\documentclass[a4paper,12pt]{article} 

\usepackage{amsmath,amssymb}
\usepackage[usenames]{color}
\usepackage{graphicx}
\usepackage{amsfonts}
\usepackage{comment}
\usepackage{pxrubrica}

\newtheorem{thm}{Theorem}[section]
\newtheorem{prop}[thm]{Proposition}
\newtheorem{lmm}[thm]{Lemma}
\newtheorem{defn}[thm]{Definition}

\newtheorem{cry}[thm]{Corollary}

\newtheorem{rmk}{Remark}[section]

\newtheorem{cmt}{Comment}
\numberwithin{equation}{section}

\newcommand{\qed}{\hfill\ensuremath{\square}\\}%

\newcommand{\pr}{\par \vspace{0mm} \noindent {\bf [Proof]} \quad}
\newcommand{\prend}{\hfill \qed}

\newcommand{\pd}{{\partial}} 
\newcommand{\gr}{{\rm gr}}

\newcommand{\wt}{{\rm wt}}

\newcommand{\Hom}{{\rm Hom}}

\newcommand{\End}{{\rm End}}

\newcommand{\Coeff}{{\rm Cf}}

\newcommand{\CC}{{\cal C}} 
\newcommand{\CD}{{\cal D}}

\newcommand{\CJ}{{\cal J}} 
 
\newcommand{\CM}{{\cal M}} 
\newcommand{\CN}{{\cal N}} 
\newcommand{\CO}{{\cal O}} 
\newcommand{\CP}{{\cal P}} 
\newcommand{\CQ}{{\cal Q}} 
 
\newcommand{\CY}{{\cal Y}}

\newcommand{\bC}{{\mathbb C}} 
\newcommand{\bD}{{\mathbb D}} 
\newcommand{\bZ}{{\mathbb Z}} 
\newcommand{\bN}{{\mathbb N}}

\newcommand{\bR}{{\mathbb R}}

\definecolor{skyblue}{rgb}{0.5,0.5,1}
\definecolor{fadegreen}{rgb}{0.9,1,0.9}

\begin{document}

\title{Associativity of fusion products of $C_1$-cofinite $\bN$-gradable modules 
of vertex operator algebra}
\author{
\begin{tabular}{c}
Masahiko Miyamoto
\footnote{Partially supported 
by the Grants-in-Aids
for Scientific Research, No.21K03195, The Ministry of Education,
Science and Culture, Japan} \cr
Institute of Mathematics, University of Tsukuba\cr
Institute of Mathematics, Academia Sinica\cr
\end{tabular}}
\date{}
\maketitle

\begin{abstract}
We prove an associative law of the fusion products $\boxtimes$ of $C_1$-cofinite $\bN$-gradable modules for a vertex operator algebra $V$. To be more precise, for $C_1$-cofinite $\bN$-gradable $V$-modules $A,B,C$ and their fusion products 
$(A\!\boxtimes\! B, \CY^{AB})$, $((A\!\boxtimes\! B)\!\boxtimes\! C, \CY^{(AB)C})$, 
$(B\!\boxtimes\! C, \CY^{BC})$, $(A\!\boxtimes\! (B\!\boxtimes\! C),\CY^{A(BC)})$ with 
logarithmic intertwining operators $\CY^{AB},\ldots,\CY^{A(BC)}$ satisfying the universal properties for 
$\bN$-gradable modules, 
we prove that four-point correlation functions 
$\langle \theta, \CY^{A(BC)}(v,x)\CY^{BC}(u,y)w\rangle$ and 
$\langle \theta', \CY^{(AB)C}(\CY^{AB}(v,x-y)u,y)w\rangle$ 
are locally normally convergent over $\{(x,y)\in \bC^2 \mid 0\!<\!|x\!-\!y|\!<\!|y|\!<\!|x|\}$. We then take their respective principal branches \vspace{-2mm}
\begin{center}
$\widetilde{F}(\langle \theta,\CY^{A(BC)}(v,x)\CY^{BC}(u,y)w\rangle)$ 
and $\widetilde{F}(\langle \theta,\CY^{(AB)C}(\CY^{AB)}(v,x-y)u,y)w\rangle)$ \vspace{-2mm}
\end{center}
on $\CD^2\!=\!\{(x,y)\in \bC^2 \mid 0\!<\!|x\!-\!y|\!<\!|y|\!<\!|x|, 
\mbox{ and } x,y,x\!-\!y\not\in \bR^{\leq 0}\}$ and then show that there is an isomorphism 
$\phi_{[AB]C}:(A\boxtimes B)\boxtimes C \to A\boxtimes (B\boxtimes C)$ such that 
\vspace{-2mm}
$$
\widetilde{F}(\langle \theta, \CY^{A(BC)}(v,x)\CY^{BC}(u,y)w\rangle)
=\widetilde{F}(\langle \phi_{[AB]C}^{\ast}(\theta), \CY^{(AB)C}(\CY^{AB}(v,x-y)u,y)w)\rangle \vspace{-2mm} $$
on $\CD^2$ 
for $\theta\in (A\boxtimes (B\boxtimes C))^{\vee}$, $v\in A$, $u\in B$, and 
$w\in C$, where $W^{\vee}$ denotes the contragredient module of $W$ and 
$\phi_{[AB]C}^{\ast}$ denotes the dual of $\phi_{[AB]C}$. 
We also prove the pentagon identity. 
\end{abstract}

\section{Introduction}
As a counterpart of a tensor product of two modules in the VOA theory, Huang and Lepowsky\cite{HL} have introduced a concept of 
fusion products $\boxtimes_{p(x)}$. 
Huang \cite{H1} has also shown that four-point correlation functions for 
three $C_1$-cofinite $\bN$-gradable modules satisfy differential equations with regular singularity and then as an application, he has proved an associative law of fusion products under some conditions, 
where a $V$-module $W$ is called to be "$C_1$-cofinite" if $\dim V/C_1(W)<\infty$ for 
$C_1(W):={\rm span}_{\bC}\{\alpha_{-1}w\in W\mid w\in W, \alpha\in V, \wt(\alpha)\geq 1\}$ and a $V$-module $W$ is called $\bN$-gradable if $W$ has a decomposition $W=\oplus_{m\in \bN}W_{(m)}$ satisfying 
$\alpha_{k}W_{(m)}\subseteq W_{(m+\wt(\alpha)-k-1)}$ for 
$\alpha\in V$ and $k\in \bZ$. 
From now on, $\CN\!\CC_1$ denotes the set of all $C_1$-cofinite $\bN$-gradable $V$-modules. 
Before we explain our results, let us show basic knowledge for $C_1$-cofinite $\bN$-gradable modules. 
If $W=\oplus_{m\in \bN}W_{(m)}\in \CN\!\CC_1$ is indecomposable, 
then as we will show in \S 2, $\dim W_{(m)}<\infty$ for all $m$ and 
there is $d=d(W)\in \bC$ such that 
$W=\oplus_{m\in \bN} W_{d+m}$ with $W_d\not=0$ and 
$\dim W_{d+m}<\infty$, where $W_s$ denotes a generalized eigenspace 
of $L(0)$ with eigenvalue $s$. 
In particular, $W^{\vee}=\oplus_{m\in \bN}\Hom(W_{(m)},\bC)
\cong \oplus_{d\in \bC}\Hom(W_{d},\bC)$. An element $w\in W_s$ is called homogeneous and denoted by $\wt(w)=s$. 
Clearly, $W_{(m)}=W_{d+m}$ defines an $\bN$-grading on $W$. 
We always choose this grading for indecomposable 
modules $A,B,C$ which come out later. We write $\gr(w)=m$ for $w\in W_{(m)}$ and call it the grade of $w$. 

For $A,B\in \CN\!\CC_1$, as the author has shown in \cite{M1}, a pair $(A\boxtimes B, \CY^{AB})$ of a $V$-module $A\boxtimes B$ 
and a logarithmic intertwining operator $\CY^{AB}\in I\binom{A\boxtimes B}{A\,\, B}$ satisfying the universal property for $\bN$-gradable modules (see \S 2) is always uniquely well-defined up to isomorphism and $A\boxtimes B$ is 
also in $\CN\!\CC_1$ and 
$\CY^{AB}$ is also a logarithmic intertwining operator.
Let $A,B,C\in \CN\!\CC_1$ and assume to be $\bN$-graded. 
Since $A\boxtimes B, B\boxtimes C\in \CN\!\CC_1$, 
we also have 
$((A\boxtimes B)\boxtimes C, \CY^{(AB)C})$ and 
$(A\boxtimes(B\boxtimes C), \CY^{A(BC)})$ and 
we denote these logarithmic intertwining operators by 
$$\begin{array}{ll}
\CY^{A(BC)}(v,x)&=\sum_{h=0}^{K_1} \CY_h^{A(BC)}(v,x)\log^h(x), \vspace{1mm}\cr
\CY^{BC}(u,y)&=\sum_{k=0}^{K_2} \CY^{BC}_k(u,y)\log^k(y),\vspace{1mm}\cr
\CY^{AB}(v,x-y)&=\sum_{h=0}^{K_3} \CY^{AB}_k(v,x-y)\log^k(x-y), \mbox{ and} \vspace{1mm}\cr 
\CY^{(AB)C}(\delta,y)&=\sum_{k=0}^{K_4} \CY_h^{(AB)C}(\delta,y)\log^h(y), 
\end{array}$$
for $v\in A, u\in B, \delta\in A\boxtimes B$ 
with formal $\bC$-power series $\CY_h^{A(BC)}(v,x)$ of $x$ with coefficients 
in $\Hom(B\boxtimes C,A\boxtimes(B\boxtimes C))$, etc.  
To simplify the notation, we use notation  
$K^2=\{ (h,k)\in \bN^2\mid h\leq K_1, k\leq K_2\}$, ${K'}^2=\{ (h,k)\in \bN^2\mid h\leq K_3, k\leq K_4\}$ and we often denote 
$A\boxtimes(B\boxtimes C)$ and $(A\boxtimes B)\boxtimes C$ by 
$A(BC)$ and $(AB)C$, respectively. 
Then for $\theta\in (A(BC))^{\vee}, \theta'\in ((AB)C)^{\vee}$, we define two four-point correlation functions and the coefficients 
of $\log^h(x)\log^k(y)$ and $\log^h(x-y)\log^k(y)$ in them as follows:
$$\begin{array}{ll}
F^{A(BC)}(\theta,v,u,w;x,y)&:=\langle \theta, \CY^{A(BC)}(v,x)\CY^{BC}(u,y)w\rangle, \vspace{1mm}\cr 
F^{(AB)C}(\theta',v,u,w;x-y,y)&:=\langle \theta', \CY^{(AB)C}(\CY^{AB}(v,x-y)u,y)w\rangle, 
\vspace{1mm} \cr
F_{h,k}^{A(BC)}(\theta,v,u,w;x,y)&:=\langle \theta, \CY_h^{A(BC)}(v,x)\CY_k^{BC}(u,y)w\rangle, \mbox{ and}\vspace{1mm} \cr
F_{h,k}^{(AB)C}(\theta',v,u,w;x-y,y)&:=\langle \theta',\CY_h^{(AB)C}(\CY_k^{AB}(v,x-y)u,y)w
\rangle.
\end{array}$$ 
Set $\Omega=(A(BC))^{\vee}\times A\times B\times C$ (and 
$\Omega'=((AB)C)^{\vee}\times A\times B\times C$). We use $\vec{\xi}$ (and $\vec{\xi}'$) to denote 
a quadruple $(\theta,v,u,w)\in 
\Omega$ (and $(\theta',v,u,w)\in \Omega'$) so that $F_{h,k}^{A(BC)}(\vec{\xi};x,y)$ denotes $F_{h,k}^{A(BC)}(\theta,v,u,w;x,y)$ and so on. 
For $\alpha\in V$, we use notation $\alpha_n^{[i]}$ to denote the action of $\alpha_n$ on the $i$-th coordinate of $\Omega$, e.g. 
$\alpha^{[2]}_n\vec{\xi}=(\theta,\alpha_nv,u,w)$. 
For $\vec{\xi}=(\theta,v,u,w)\in \Omega$, 
define $\gr^{234}(\vec{\xi})=\gr(v)+\gr(u)+\gr(w)$ and the total grade 
$\gr(\vec{\xi})=\gr(\theta)+\gr^{234}(\vec{\xi})$.

The main purpose in this paper is to show the associative law of fusion products (Theorem 4.1) 
for $A,B,C\in \CN\!\CC_1$. We may assume that $A,B,C$ are indecomposable. 
As a proof, we will construct a surjective homomorphism 
$\phi_{[AB]C}:(A\boxtimes B)
\boxtimes C\to A\boxtimes(B\boxtimes C)$ by starting from 
$F^{A(BC)}(\vec{\xi};x,y)$. The reverse homomorphism is given by 
starting from $F^{(AB)C}(\vec{\xi}';x,y)$ with the similar argument. 
Our proofs are based on Huang's ideas in \cite{H1}. The differences from \cite{H1} are that we will treat formal $\bC$-power series $F_{h,k}^{A(BC)}(\vec{\xi};x,y)$ and 
show that their modified versions satisfy two Fuchsian systems on discs 
by fixing one variable $y=y_0\not\in \bR^{\leq 0}$ or $x=x_0\not\in \bR^{\leq 0}$. 
Furthermore, in order to expand $F^{A(BC)}(\vec{\xi};x,y_0)$ in a $\bC$-formal power series of 
$(x-y_0)$ with logarithm functions $\log^h(x-y_0)$ as a component of a solution of 
Fuchsian system,   
we restrict the variable $x$ into a domain $\CD^2_{(x,y_0)}\!=\!\{x\in \bC
 \mid 0\!<\!|x\!-\!y_0|\!<\!|y_0|\!<\!|x|, \mbox{ and }x,y_0,x\!-\!y_0\not\in \bR^{\leq 0} \}$ and take a principal branch $\widetilde{F}^{A(BC)}(\vec{\xi};x,y_0)$ of $F^{A(BC)}(\vec{\xi};x,y_0)$. 
A key result we get from our Fuchsian system is that 
there is a finite set $\Delta\subseteq \bC$ such that the powers of $(x-y_0)$ in 
these expansions are all in $\Delta-\gr(v)-\gr(u)+\bN$ for any $\vec{\xi}$, 
which supports an existence of an $\bN$-gradable $V$-module isomorphic to $A\boxtimes B$.

In order to get Fuchsian systems, we consider finite spaces. 
Let $P_{W}$ denote a 
complement of $C_1(W)$ in $W$ and let $\tilde{P}_W$ be a finite dimensional 
subspace of $W$ containing $P_W$ and spanned by homogeneous elements.  
For $N\in \bN$, we define $W^{\vee}_{(\leq N)}=\oplus_{m=0}^{N}\Hom(W_{(m)},\bC)$ and 
set $\Omega_N=\{(\theta,v,u,w)\in (A(BC))^{\vee}_{(\leq N)}\times A\times B\times C\}$. 
Choose bases $J_{\tilde{P}_A}=\{v^i\mid i\in \CP_A\}, J_{\tilde{P}_B}=\{u^j\mid j\in \CP_B\}, 
J_{\tilde{P}_C}
=\{w^q\mid q\in \CP_C\}, J_N=\{\theta^p\mid p\in\CP_N \}$, of $\tilde{P}_A, \tilde{P}_B, \tilde{P}_C, 
(A(BC))^{\vee}_{(\leq N)}$ consisting of homogeneous elements and define  
$\CJ_{N,\tilde{P}_A,\tilde{P}_B,\tilde{P}_C}
=J_{N}\times J_{\tilde{P}_A}\times J_{\tilde{P}_B}\times J_{\tilde{P}_C}$, which is a finite set. 
We simply denote it by $\CJ_N$ and denote $\CJ_{N,P_A,P_B,P_C}$ by $\CJ_N^0$. 
We always assume $\CJ^0_N\subseteq \CJ_N$. 
We similarly define 
$\Omega'_{N}$ and finite sets $\CJ'_N$ and ${\CJ'}^0_{N}$ for $((AB)C)^{\vee}$. 
Set 
$$\begin{array}{ll}
G^{A(BC):y}_{h,k}(\vec{\xi};x,y)
&:=\langle \theta,\CY_h^{A(BC)}(v,x)\CY_k^{BC}(u,y)w\rangle y^{\gr^{234}(\vec{\xi})},\vspace{1mm}\cr
G^{A(BC):x-y}_{h,k}(\vec{\xi};x,y)&:=
\langle \theta, \CY_h^{A(BC)}(v,x)\CY_k^{BC}(u,y)w\rangle (x-y)^{\gr^{234}(\vec{\xi})}, \vspace{1mm}\cr
G^{(AB)C:y}_{h,k}(\vec{\xi}';x,y)
&:=\langle \theta',\CY_h^{(AB)C}(\CY_k^{AB}(v,x-y)u,y)w\rangle y^{\gr^{234}(\vec{\xi})}, 
\mbox{ and}\vspace{1mm}\cr
G^{(AB)C:x-y}_{h,k}(\vec{\xi}';x,y)&:=
\langle \theta', \CY_h^{(AB)C}(\CY_k^{AB}(v,x-y)u,y)w\rangle (x-y)^{\gr^{234}(\vec{\xi})}.  
\end{array}$$  

Then we will obtain the following reduction theorem. \\

\noindent 
{\bf Theorem 3.1 for $G^{A(BC)}$} \quad {\it 
For $\vec{\xi}=(\theta,v,u,w)\in \Omega_N$, $\alpha\in V$, and $x_0\not=0\not=y_0$, \\
(1) $G_{h,k}^{A(BC):y}(\vec{\xi};x_0,y)$ is a linear combination of 
$\{G_{h,k}^{A(BC):y}(\vec{\mu};x_0,y)\mid \vec{\mu}\in 
\CJ^0_N\}$ 
with coefficients 
in $\bC[\iota_{x_0,y}\{(x_0-y)^{-1}\}] [y]\subseteq \bC[[y]]$ and \\
(2) $G_{h,k}^{A(BC):x-y}(\vec{\xi};x,y_0)$ 
is a linear combination of 
$\{G_{h,k}^{A(BC):x-y}(\vec{\mu};x,y_0)\mid 
\vec{\mu}\in \CJ^0_N\}$ with coefficients 
in $\bC[x,\iota_{y_0,x-y_0}\{x^{-1}\}][x-y_0]\subseteq \bC[[x-y_0]]$. \\
Furthermore, for the residue classes of coefficients modulo $\bC[[x-y_0]](x-y_0)$, we have:
$$\begin{array}{ll}
G_{h,k}^{A(BC):x-y}(\alpha_{-1}^{[2]}\vec{\xi};x,y_0)\equiv \sum_{j=0}^{\infty} 
G_{h,k}^{A(BC):x-y}(\alpha_j^{[3]}\vec{\xi};x,y_0)
&\pmod{\bC[[x\!-\!y_0]](x\!-\!y_0)},\cr
G_{h,k}^{A(BC):x-y}(\alpha_{-1}^{[3]}\vec{\xi};x,y_0)\equiv \sum_{j=0}^{\infty}(-1)^{j+1} G_{h,k}^{A(BC):x-y}( \alpha_j^{[2]}\vec{\xi};x,y_0) \hspace{-6mm}
&\pmod{\bC[[x\!-\!y_0]](x\!-\!y_0)},\cr
G_{h,k}^{A(BC):x-y}(\alpha_{-1}^{[4]}\vec{\xi};x,y_0)\equiv 0 &\pmod{\bC[[x\!-\!y_0]](x\!-\!y_0)}. 
\end{array}$$
}
\vspace{1mm}

Using $L(-1)$-derivative properties: 
$F^{A(BC)}(L(-1)^{[3]}\vec{\xi};x,y)=\frac{\pd}{\pd y}F^{A(BC)}(\vec{\xi};x,y)$ 
and 
$F^{A(BC)}(L(-1)^{[2]}\vec{\xi};x,y)=\frac{\pd}{\pd x}F^{A(BC)}(\vec{\xi};x,y)$, 
we will show that there are 
$\lambda^{34,\vec{\xi},h,k}_{\vec{\mu},p,q}(x_0,y)\in \bC[\iota_{x_0,y}\{(x_0-y)^{-1}\}][y]$ and $\lambda^{23,\vec{\xi},h,k}_{\vec{\mu},p,q}(x,y_0)\in 
\bC[x,\iota_{y_0,x-y_0}\{x^{-1}\}][x-y_0]$ 
such that 
$$\begin{array}{l}
 \frac{\pd}{\pd y}G^{A(BC):y}_{h,k}(\vec{\xi};x_0,y)
=\frac{1}{y}\sum_{\vec{\mu}\in \CJ^0_N}\sum_{(p,q)\in K^2}\lambda^{34,\vec{\xi},h,k}_{\vec{\mu},p,q}(x_0,y)G^{A(BC):y}_{p,q}(\vec{\mu};x_0,y) \mbox{ and}\cr
\frac{\pd}{\pd x}G^{A(BC):x-y}_{h,k}(\vec{\xi};x,y_0)
=\frac{1}{x-y_0}\sum_{\vec{\mu}\in \CJ^0_N}\sum_{(p,q)\in K^2}
\lambda^{23,\vec{\xi},h,k}_{\vec{\mu},p,q}(x,y_0)G^{A(BC):x-y}_{p,q}(\vec{\mu};x,y_0) 
\end{array}$$
as formal $\bC$-power series. Set $r\!=\!|\CJ_N||K^2|$ and define matrix-valued 
functions \\$\Lambda^{23}(x,y_0)\!=\!\left(\lambda^{23,\vec{\xi},h,k}_{\vec{\mu},p,q}(x,y_0)\right)_{\!(\vec{\xi},h,k),(\vec{\mu},p,q)\in \CJ_N\times K^2}
\!\!\in\! M_{r\times r}(\CO(\bD_{|y_0|}(y_0)))$ and \\
$\Lambda^{34}(x_0,y)\!=\!
\left(\lambda^{34,\vec{\xi},h,k}_{\vec{\mu},p,q}(x_0,y)\right)_{(\vec{\xi},h,k),(\vec{\mu},p,q)\!\!\in \CJ_N\times K^2}\!\in\! M_{r\times r}(\CO(\bD_{|x_0|}(0)))$. Here and hereafter 
$\bD_R(z_0)=\{z\in \bC \mid |z\!-\!z_0|\!<\!R\}$ and $\CO(U)$ denotes the set of holomorphic functions on $U$. 
We note $\lambda^{23,\vec{\xi},h,k}_{\vec{\mu},p,q}\!=\!0\!=\!\lambda^{34,\vec{\xi},h,k}_{\vec{\mu},p,q}$ if $\vec{\mu}\not\in \CJ^0_N$. 
Then we get: \\

\noindent
{\bf Theorem 3.2} [Differential systems] 
\quad{\it
Fix $x_0\not=0$ and $N\in \bN$. Then 
a vector valued function $G^y(x_0,y):
=\left(G^{A(BC):y}_{h,k}(\vec{\xi};x_0,y)\right)_{ \vec{\xi}\in \CJ_{N},(h,k)\in K^2}$
satisfies a Fuchsian system 
$$\frac{d}{dy}G^y(x_0,y)=\frac{\Lambda^{34}(x_0,y)}{y}G^y(x_0,y)\eqno{(D1)}$$ 
over $\bD_{|x_0|}(0)$ with polar locus $\{0\}$. 
Similarly, for a fixed $0\not=y_0\in \bC$ and $N\in \bN$, 
a vector valued function 
$G^{x-y}(x,y_0):=\left( G^{A(BC):x-y}_{h,k}(\vec{\xi};x,y_0)\right)_{\vec{\xi}\in \CJ_{N}, (h,k)\in K^2}$ satisfies a Fuchsian system 
$$\frac{d}{dx}G^{x-y}(x,y_0)=\frac{\Lambda^{23}(x,y_0)}{x-y_0}G^{x-y}(x,y_0) \eqno{(D2)}$$ 
over $\bD_{|y_0|}(y_0)$ with polar locus $\{y_0\}$. 
Furthermore, the sets of nonzero eigenvalues (without multiplicities) 
of the constant matrices $\Lambda^{34}(x_0,0)$ of $(D1)$ and 
$\Lambda^{23}(y_0,y_0)$ of $(D2)$ are determined by 
only the choices of the bases $J_{P_A}$ and $J_{P_B}$ of $P_A$ and $P_B$. 

We also have similar results for 
$\bar{G}^{y}(x_0,y)=\left( G^{(AB)C:y}_{h,k}(\vec{\xi}';x_0,y) \right)_{\vec{\xi}'\in \CJ'_N,(h,k)\in {K'}^2}$ and 
$\bar{G}^{x-y}(x,y_0)=\left(G^{(AB)C:x-y}_{h,k}(\vec{\xi}',x,y_0)\right)_{\vec{\xi}'\in 
\CJ'_N,(h,k)\in {K'}^2}$. }\\
\vspace{2mm}

As a corollary of $(D1)$ for $G^{A(BC)}$ and $(D2)$ for $G^{(AB)C}$, we will obtain the following: \\

\noindent
{\bf Corollary 3.3}\quad {\it 
$F_{h,k}^{A(BC)}(\vec{\xi};x,y)$ and $F_{h,k}^{(AB)C}(\vec{\xi};x-y,y)$ 
are locally normally convergent on $\{(x,y)\in \bC^2 \mid 0<|y|<|x|\}$ and 
$\{(x,y)\in \bC^2 \mid 0<|x-y|<|y|\}$, respectively.}\\

Set $\CD^2=\{(x,y)\in \bC^2 \mid 0<|x\!-\!y|<|y|<|x| \mbox{ and }x,y,x\!-\!y\not\in \bR^{\leq 0}\}$. 
Let $y_0\not\in \bR^{\leq 0}$ and we take a principal branch $\widetilde{G}^{A(BC):x-y}_{h,k}(\vec{\xi};x,y_0)$
of $G^{A(BC):x-y}_{h,k}(\vec{\xi};x,y_0)$ on a small domain 
$\CD^2_{(x,y_0)}\!=\!\{x\in \bC\mid (x,y_0)\in \CD^2\}$. 
As we will explain in Proposition \ref{Sol} in \S 2, there is a finite subset 
$\Delta'=\{d_1,...,d_p\} \subseteq \bC$ which depends only on the choice of $J_{P_A}$ and $J_{P_B}$ such that all components of solutions of 
the Fuchsian system (D2) with polar locus $\{y_0\}$ are written as  
$$\sum_{d\in \Delta'}\sum_{m\in \bN}\sum_{t=0}^K r_{d+m,t}(x-y_0)^{d+m}\log^t(x-y_0) $$ 
with $r_{d+m,t}\in \bC$ 
at $\bD_{|y_0|}(y_0)=\{x\in \bC \mid 0<|x-y_0|<|y_0|\}$. Since this is true for every $y=y_0\not\in \bR^{\leq 0}$ and $\vec{\xi}$, there is $r^{h,k}_{d+m,t}(\vec{\xi},y)\in \bC$ and 
$K(\vec{\xi})\in \bN$ such that 
$$\widetilde{G}^{A(BC):x-y}_{h,k}(\vec{\xi};x,y)=
\sum_{d\in \Delta'}\sum_{m\in\bN}\sum_{t=0}^{K(\vec{\xi})} r^{h,k}_{d+m,t}(\vec{\xi},y)(x-y)^{d+m}
\log^t(x-y). $$
Multiplying it by $(x-y)^{-\gr^{234}(\vec{\xi})}\log^h(x)\log^k(y)$ for $h,k$ and 
taking a sum of them and replacing $\log(x)$ by 
$\sum_{j=0}^{\infty} \frac{(-1)^j}{j+1}(\frac{x-y}{y})^{j+1}+\log(y)$, 
we obtain an expansion of a principal branch 
$\widetilde{F}(\langle \theta, 
\CY^{A(BC)}(v,x)\CY^{BC}(u,y)w\rangle)$ of $F^{A(BC)}(\vec{\xi};x,y)$ on 
$\CD^2$. 
We denote it by $\widetilde{F}^{A(BC)}(\vec{\xi};x,y)$. 
Since $\dim P_C\!<\!\infty$, setting $\Delta\!=\Delta'!-\!{\rm max}\{\gr(w)\mid w\in\! P_C\}$, we can write  
$$
\widetilde{F}^{A(BC)}(\vec{\xi};x,y)=\sum_{t=0}^{K(\vec{\xi})} \sum_{-s-1\in \Delta+\bN-\gr(v)-\gr(u)} g_{s,t}(\vec{\xi};y)(x-y)^{-s-1}\log^t(x-y)\eqno{(4.5)}$$
for $\vec{\xi}=(\theta,v,u,w)\in \Omega$ with $w\in P_C$ on $\CD^2$.  
Although our argument is depending on the choice $\CJ_N$, 
we will show that $(4.5)$ holds for every $\vec{\xi}\in \Omega$ with the same 
$\Delta$.

Since $g_{s,t}(\vec{\xi};y)$ is multi-linear on $\vec{\xi}\in (A(BC))^{\vee}\times 
A\times B\times C$, we will view it as a map $A\times B 
\to \Hom(C, A(BC))\otimes \CO(\bC\backslash \bR^{\leq 0})$, 
in other words, using formal operators $\CY^{3}$ and $\CY^{4}$, 
we can rewrite RHS of (4.5) into $\langle \theta, \CY^{3}(\CY^{4}(v,x-y)u,y)w\rangle$ such that 
$$ \widetilde{F}^{A(BC)}(\vec{\xi};x,y)=\langle \theta, \CY^{3}(\CY^{4}(v,x-y)u,y)w\rangle \eqno{(4.9)}$$
on $\CD^2$ for $\vec{\xi}\in \Omega$. 
In \S 4, using these notation, we will prove that $\CY^3$ and $\CY^4$ are logarithmic surjective intertwining 
operators, which implies the existence of surjective homomorphism 
$(A\boxtimes B)\boxtimes C\to A\boxtimes(B\boxtimes C)$ and then we will get 
the main theorem (Theorem 4.1). \\

We will also prove another axiom of monoidal category called pentagon axiom in \S 5. \\

\noindent
{\bf Acknowledgments.} The author wishes thank Toshiyuki Abe, Scott 
Carnahan, Lam Ching Hung, Hiroshi Yamauchi and Atsushi Matsuo for their 
useful advice and comments. He also 
wishes to thank Mathematical Institute of Academia Sinica for hospitality. 

\section{Notation and preliminary results}

\subsection{Differential systems}
We quote the following statements from \cite{Heu}.
For the proofs, see \cite{Heu}. See also \cite{IY}. 

\begin{defn} A meromorphic system of linear differential equations is a 
differential equation of the form 
$\frac{d}{dz}Y=AY$, 
where $A\in M_{r\times r}(\CM(U))$ is a square matrix whose coefficients 
are meromorphic functions over some complex domain 
$U\subseteq \bC$ (open, nonempty). 
\end{defn}

\noindent 
$\bullet$ A singularity of such a meromorphic system is a point in $U$ corresponding to a pole of 
(at least one of the coefficients of) $A$. \\
$\bullet$ For technical reasons, one usually fixes a discrete subset $\Sigma\subseteq U$, and 
speaks of meromorphic systems $\frac{d}{dz}Y=AY$ over $U$ with polar locus $\Sigma$ if 
the "true polar locus" (i.e., the set of poles of $A$) is contained in $\Sigma$. \\
$\bullet$ A meromorphic system $\frac{d}{dz}Y=\frac{A}{z-z_0}Y$ with a holomorphic matrix function $A$ of size 
$r\times r$ defined on $U$ is called a Fuchsian system and its singularity at $z=z_0$ is called a Fuchsian singularity.

\begin{defn} Let $\frac{d}{dz}Y=AY$ be a meromorphic system of rank $r$ over $U\subseteq \bC$. 
Let $U'\subseteq U$ be open. A fundamental solution of $\frac{d}{dz}Y=AY$ over $U'$ is a matrix-valued meromorphic function 
$\Phi\in M_{r\times r}(\CM(U'))$ such that 
$$\frac{d}{dz}\Phi(z)=A(z)\Phi(z) \quad  \forall z\in U', $$
and moreover 
$$  \det(\Phi)(z)\not=0 \quad \forall z\in U'. $$
\end{defn}

From now on, let $z_0\in U$ and 
let $A=\sum_{k\geq 0}(z-z_0)^kA_k$ with $A_k\in M_{r\times r}(\bC)$ 
and we will consider 
$$ \frac{d}{dz}Y=\frac{A}{z-z_0}Y,  \eqno{(2.1)}$$ 
which we will call Fuchsian system. Set $\bD_R(z_0)=\{z\in \bC 
\mid |z-z_0|<R\}$. 
If $A=A_0\in M_{r\times r}(\bC)$ and $z_0\in U$, 
then 
$$\frac{d}{dz}Y=\frac{A_0}{z-z_0}Y \eqno{(2.2)}$$ 
is called Euler system. In this case, 
there is $G\in GL_r(\bC)$ such that $GA_0G^{-1}$ has a 
Jordan standard form and 
$\frac{d}{dz}(GY)=G\frac{d}{dz}Y=\frac{GA_0G^{-1}}{z-z_0}(GY)$. 

\begin{lmm}\label{Jordan} 
If an Euler system (2.2) has a Jordan cell $A_0=J_r(a)$, 
that is, $Y={}^t(f_1,\ldots,f_r)$ satisfies $(z-z_0)\frac{d}{dz}f_i=af_i+f_{i+1}$ for $i=1,...,r$ with 
$f_{r+1}=0$, then 
$\Phi=\{(z-z_0)^a\times {}^t(\log^s(z-z_0),...,\log(z-z_0),1,0,\cdots,0)\mid s=0,...,r-1\}$ 
is a fundamental solution, where ${}^t(\cdots)$ denotes the transpose of $(\cdots)$. 
\end{lmm}

Therefore, in generally, all components of solutions of Euler system (2.2) are written as 
$$f(z)=\sum_{i=1}^r\sum_{t=0}^{r}r_{i,t} (z-z_0)^{d_i}\log^t(z-z_0) 
\eqno{(2.3)}$$
with $r_{i,t}\in \bC$ and eigenvalues $d_i$ of $A_0$. 
We also note that if $\frac{d}{dz}Y=\frac{A}{z-z_0}Y$, then 
$$\frac{d}{dz}((z-z_0)^dY)=\frac{(A+dE_r)}{z-z_0}((z-z_0)^dY),  \eqno{(2.4)} $$
where $E_r$ is the identity matrix of size $r$. \\

\begin{defn}
Let $\frac{d}{dz}Y=AY$ be a meromorphic system over $U$ with polar locus $\Sigma$. \\
$\bullet$ A holomorphic gauge transformation over $U'\subseteq U$ is the change of variable 
$Y=\Delta Z$, 
yielding a meromorphic system $\frac{d}{dz}Z=BZ$ over $U'$ with polar locus 
$\Sigma\cap U'$, where 
$\Delta\in GL_r(\CO(U'))$, 
i.e. $\Delta$ is a holomorphic matrix function over $U'$ with non-vanishing determinant.\\
$\bullet$ A meromorphic gauge transformation 
over $U'\subseteq U$ is the change of variable $Y=\Delta Z$, 
yielding a meromorphic system $\frac{d}{dz}Z=BZ$ over $U'$ with polar locus 
$Z\cap U'$, where 
$\Delta\in GL_r(\CM(U'))$, 
i.e. $\Delta$ is a meromorphic matrix function over $U'$ whose determinant is not 
identically zero on any connected component of $U'$, such that moreover, 
$\Delta|_{U'\backslash \Sigma}\in M_{r\times r}(\CO(U'\backslash \Sigma))$ 
and 
$\Delta^{-1}|_{U'\backslash \Sigma}\in M_{r\times r}(\CO(U'\backslash \Sigma))$. 
\end{defn}

\noindent
$\bullet$ As usual, two meromorphic systems over $U$ with polar locus $\Sigma$ are said to be 
holomorphically (resp. meromorphically) gauge equivalent, if they are related via a global holomorphic (resp. meromorphic) gauge transformation. \\

\begin{defn}[Definition 3.6 in \cite{Heu}]
A Fuchsian system $\frac{d}{dz}Y=\frac{B}{z}Y$ of rank $r$ over $\bD$, with polar locus $\{0\}$, is said to be of Levelt normal form if 
$B=\sum_{k\geq 0}B_kz^k$ such that for any $k\in \bN$, the matrix $B_k\in M_{r,r}(\bC)$ satisfies 
${\rm ad}(B_{0,s})(B_k)=kB_k$, where $B_0=B_{0,s}+B_{0,n}$ is the Dunford decomposition of $B_0$ into semisimple part $B_{0,s}$ and nilpotent part $B_{0,n}$ with $B_{0,s}B_{0,n}=B_{0,n}B_{0,s}$. 
\end{defn}

\begin{lmm} 
The eigenvalues of $B(1)=\sum_{k=0}^{\infty} B_k$ are equal to the eigenvalues of $B_0$.  
\end{lmm}

\pr 
We may assume that $B_{0,s}$ is a diagonal matrix $(\lambda_1,...,\lambda_r)$ and 
$B_{0,n}$ is an upper triangular matrix. Furthermore, we may assume $\Re(\lambda_1)\geq\cdots\geq \Re(\lambda_r)$. Then since ${\rm ad}(B_{0,s})E_{i,j}=(\lambda_i-\lambda_j)E_{i,j}$, 
${\rm ad}(B_{0,s})B_k=kB_k$ means that $B_k$ are strictly upper triangular and  
$B(1)=\sum_{k=0}^{\infty} B_k$ is an upper triangular matrix with diagonal entries 
$\{\lambda_1,\ldots,\lambda_r\}$. \qed

\begin{thm}[Theorem 3.8 in \cite{Heu}] 
Let $\frac{d}{dz}Y=\frac{A}{z}Y$ be a Fuchsian system over 
$\bD$. This system is holomorphically 
gauge equivalent to some Fuchsian system 
$\frac{d}{dz}\widehat{Y}=\frac{B}{z}\widehat{Y}$ of Levelt normal form over $\bD$.
\end{thm}

\begin{rmk}
Under the setting of Theorem 2.7, 
the set of eigenvalues of $B(0)$ coincides with 
the set of eigenvalues of $A(0)$ counted with multiplicity.  
\end{rmk}

\begin{thm}[Lemma 3.7 in \cite{Heu}]
Let $B_{0,s}\in M_{r\times r}(\bC)$ be semisimple. Then 
there exists the unique semisimple matrix $L$ such that 
$L$ commutes with $B_{0,s}$ and has only integer eigenvalues and 
any eigenvalues $\mu$ of $B_{0,s}-L$ satisfies $\Re(\mu)\in [0,1[$. 
Moreover, if $\frac{d}{dz}Y=\frac{B}{z}Y$ is a Fuchsian system of Level normal form with 
$B(0)=B_{0,s}+B_{0,n}$ as above, then $B=z^LB(1) z^{-L}$.  
\end{thm}

\begin{rmk}
As mentioned in the proof of Corollary 3.9 in \cite{Heu}, 
under the assumption in Theorem 2.8, 
$\tilde{Y}:=z^{-L} Y$ satisfies an Euler system 
$\frac{d}{dz}\tilde{Y}=\frac{B(1)-L}{z}\tilde{Y}$. 
\end{rmk}

Using the these results, we get the following result. 

\begin{prop}\label{Sol}
Let $\frac{d}{dz}Z=\frac{B}{z}Z$ be a Fuchsian system with $B\in M_{r\times r}(\CO(U))$ and let 
$\{d_1,\cdots,d_r\}$ be the set of eigenvalues of the constant matrix $B_0$. 
Choose $\ell_i\in \bZ$ so that $0\leq d_i-\ell_i<1$.  
Then all components of solutions of $\frac{d}{dz}Z=\frac{B}{z}Z$ 
are written as 
$$
\sum_{j=1}^r\sum_{i=1}^r\sum_{m\in \bN}\sum_{t=0}^r r^j_{i,m,t} z^{d_j-\ell_j+\ell_i+m}\log^t(z) 
\eqno{(2.5)}$$ 
with $r^j_{i,m,t}\in \bC$.
\end{prop}

\pr 
If $\frac{d}{dz}Z=\frac{A}{z}Z$ and $\frac{d}{dz}Y=\frac{B}{z}Y$ are holomorphically equivalent, then 
all components of solutions of one are linear combinations of components of solutions of the other with coefficients in 
$\bC[[z]]$ and so if all components of solutions of $\frac{d}{dz}Z=\frac{A}{z}Z$ 
satisfy (2.5), then those of solutions of 
$\frac{d}{dz}Y=\frac{B}{z}Y$ also satisfy (2.5). 
Therefore, we may assume that 
$\frac{d}{dz}Z=\frac{B}{z}Z$ is a Fuchsian system of Levelt normal form over $\bD$ by Theorem 2.7. 
Furthermore we may assume that 
$B_{0,s}={\rm diag}(d_1,...,d_r)$ satisfying 
$\Re(d_1)\geq \cdots \geq \Re(d_r)$ and $B_{0,n}$ is an upper triangular matrix. In this case,  
$L$ in Theorem 2.8 is a diagonal matrix ${\rm diag}(\ell_1,...,\ell_r)$. Then by Theorem 2.8 and Remark 2.1, we have 
$B=z^LB(1)z^{-L}$ and $\tilde{Y}:=z^{-L} Y$ satisfies an Euler system 
$\frac{d}{dz}\tilde{Y}=\frac{B(1)-L}{z}\tilde{Y}$. 
By Lemma 2.6, $B(1)$ is an upper triangular matrix with 
diagonal entries $\{d_1,...,d_r\}$ and 
the eigenvalues of $B(1)-L$ are $d_1-\ell_1,\ldots,d_r-\ell_r$. 
Hence all components of a solution $\tilde{Y}$ 
are written as 
$\sum_{j=1}^r\sum_{t\in \bN} r^j_{t} z^{d_j-\ell_j}\log^t(z)$ 
by Lemma 2.3. In this case, 
since $Y=z^L\tilde{Y}$, all components of solutions of 
$\frac{d}{dz}Y=\frac{B}{z}Y$ are written as  
$\sum_{i,j=1}^r\sum_{t\in \bN} r^j_{i,t}z^{d_j-\ell_j+\ell_i}\log^t(z)$, 
as we desired.  
\prend

Regarding the radius of convergence of such a solution, we quote the following theorem. 

\begin{thm}[Theorem 1.6 \cite{Heu}]\label{Radius}
Let $A=\sum_{k\in \bN}A_k(z-z_0)^k$ with $A_k\in M_{r\times r}(\bC)$ be a matrix power series with radius of convergence $R>0$. Let 
$Y=\sum_{k=0}^{\infty} Y_k(z-z_0)^k$ with $Y_k\in M_{r\times r}(\bC)$ be a formal power series satisfying 
$\frac{d}{dz}Y=\frac{A}{z-z_0}Y$. Then the radius of convergence of $Y$ is at least $R$.
\end{thm}

\subsection{Intertwining operators among $C_1$-cofinite modules}
Let $A,B,C$ be $V$-modules. Let $I\binom{C}{A\,\,B}$ denote the set of logarithmic intertwining 
operators of type $\binom{C}{A\,\,B}$. We introduce a generalized concept of (logarithmic) 
intertwining operators. For a nonempty domain $U\subseteq \bC\backslash\{0\}$, we define 
a local intertwining operator $\CY$ of type $\binom{C}{A\,\, B}$ on $U$ as  
a linear map $A\otimes B\to C\otimes \CO(U)$ , 
we denote it by $\langle\theta, \CY(v,z)u\rangle$ for $v\in A, u\in B, 
\theta\in C^{\vee}$ and $z\in U$, satisfying \\
${\rm (I\,\, 1)}$ [commutativity]\quad 
$\langle \theta, \alpha_m\CY(v,z)u-\CY(v,z)\alpha_mu\rangle
=\sum_{j=0}^{\infty}\binom{m}{j}\langle \theta, \CY(\alpha_jv,z)u\rangle z^{m-j}$,\\
${\rm (I\,\, 2)}$ [associativity]\quad \\
\mbox{}\qquad $\langle \theta, \CY(\alpha_mv,z)u\rangle
=\sum_{j=0}^{\infty}\binom{m}{j}(-1)^j\langle \theta, \{\alpha_{m-j}z^j\CY(v,z)u
-(-1)^m \CY(v,z) \alpha_j z^{m-j}u\}\rangle$, and  \\
${\rm (I\,\, 3)}$ [$L(-1)$-derivative property] $\langle \theta, \CY(L(-1)v,z)u\rangle
=\frac{d}{dz}\langle \theta,\CY(v,z)u\rangle$\\
for $\alpha\in V$ and $m\in \bZ$, 
where $\langle \theta, \alpha_m\CY(v,z)u\rangle$ denotes $\langle (\alpha_m)^{\ast}\theta, \CY(v,z)u\rangle$ with a contragredient operator $(\alpha_m)^{\ast}$ of $\alpha_m$. \\

An aim in this subsection is to show that every local intertwining operator for $A,B,C\in \CN\CC_1$ on a domain $U$ is 
a branch of a logarithmic intertwining operator on $U$.

Let $W\in \CN\!\CC_1$. Since $C_1(W)={\rm Span}_{\bC}\{\alpha_{-1}(P_W+C_1(W))\mid \alpha\in V_{\geq 1}\}$, we have: 

\begin{lmm}
$W={\rm Span}_{\bC}\{ 
\alpha^1_{-1}\cdots \alpha^k_{-1}w \mid w\in P_W, \alpha^j\in V_{\geq 1}, k\geq 0\}$.  
In particular, $W$ is a finitely generated $V$-module. 
\end{lmm}

Therefore, there are finitely many $d_1,...,d_t\in \bC$ such that 
$W=\oplus_{j=1}^t  \oplus_{n\in \bN} W_{d_j+n}$ and $\dim W_s<\infty$ for $s\in \bC$, where $W_s$ denotes a generalized eigenspace in $W$ of $L(0)$ with eigenvalue 
$s$. In particular, if $W$ is indecomposable, then $W=\oplus_{m\in \bN}W_{d+m}$ for some 
$d\in \bC$ and $W_d\not=0$. We use the notation 
$\wt(w)=s$ for $w\in W_s$ and $\wt(\alpha_n)=\wt(\alpha)-n-1$ 
for $\alpha\in V$. Define $\wt()1\in \End_{\bC}(W)$ 
by $(\wt()1)w=\wt(w)w$ for $w\in W$.  
Since $W_{(m)}\subseteq {\rm Span}_{\bC}
\{ \alpha^1_{-1}\cdots \alpha^k_{-1}w 
\mid w\in P_{W}, \alpha^j\in V_{\geq 1}, 
\sum_{j=1}^k \wt(\alpha^j)\leq m\}$, we obtain: 

\begin{lmm}\label{Upper} 
For $m\in \bN$, there is $\mu(m)\in \bN$ such that $\dim W_{(m)}\leq \mu(m)\times \dim P_W$. Moreover, the restricted dual $W^{\vee}=\oplus_{m\in \bN} \Hom(W_{(m)},\bC)$ of an $\bN$-graded module $W=\oplus W_{(m)}\in\CN\CC_1$ 
is isomorphic to the contragredient module $\oplus_{s\in \bC}\Hom(W_s,\bC)$ of 
$W=\oplus_{s\in \bC}W_s$. In particular, $W^{\vee}$ does not depend on the choice 
of an $\bN$-grading on $W$. 
\end{lmm}

Since $\wt(\alpha_nw)=\wt(w)+\wt(\alpha)-n-1$ and 
$L(0)\alpha_nw=(\wt(\alpha)-n-1)\alpha_nw+\alpha_nL(0)w
=(\wt(\alpha_nw)\alpha_nw+\alpha_n(L(0)-\wt(w))w$, 
we get:

\begin{lmm}\label{Nil}
$L(0)-\wt()1$ is a $V$-homomorphism. 
If $W\in \CN\!\CC_1$, then there is $K\in \bN$ such that $(L(0)-\wt()1)^KW=0$ 
and $(L(0)-\wt()1)^KW^{\vee}=0$. 
\end{lmm}

\pr 
As we have just proved, $L(0)\!-\!\wt()1$ is a $V$-homomorphism.  
Since $\dim P_W\!<\!\infty$, there is $n\!\in\! \bN$ such that 
$P_W\!\subseteq \!\oplus_{i=0}^{n} W_{(i)}$, which is $L(0)\!-\!\wt()1$-invariant and 
$\dim \oplus_{i=0}^nW_{(i)}\!<\!\infty$. Hence 
$(L(0)\!-\!\wt()1)^K(\oplus_{i=0}^nW_{(i)})=0$ for some $K$. 
Since $W$ is generated from $P_W$ by the actions of $V$, 
we have $(L(0)-\wt()1)^KW=0$. From $\langle (L(0)-\wt()1)^KW^{\vee}, W\rangle
=\langle W^{\vee}, (L(0)\!-\!\wt()1)^KW\rangle=0$, we also 
obtain $(L(0)\!-\!\wt()1)^KW^{\vee}=0$.  \qed

\begin{lmm}
Let $A,B,C\in \CN\CC_1$ and let 
$\CY$ be a local intertwining operator of type $\binom{C}{A\,\, B}$ on a domain $U$ 
with a form $\langle \theta, \CY(v,z)u\rangle=\sum_{h=0}^{K(\theta,v,u)}\langle \theta, \CY_h(v,z)u\rangle \log^h(z)$ with $\bC$-formal power series $\CY_h(v,z)u$ and $K(\theta,v,u)\in \bN$. Then there is $K\in \bN$ such that 
$K(\theta,v,u)\leq K$.  
\end{lmm}

\pr
Since 
$$\sum_{h=0}^{K(\theta,v,u)} \langle \theta, \CY_h(L(-1)v,z)u\rangle \log^h(z)
=\frac{\pd}{\pd z}\{\sum_{h=0}^{K(\theta,v,u)} \langle \theta, \CY_h(v,z)u\rangle \log^h(z)\},$$ 
we have $\CY_{k+1}(v,z)u=\CY_k(L(-1)v,z)u-\frac{\pd}{\pd z}(\CY_k(v,z)u)\times \frac{z}{k+1}$. 
We also get  
$$\begin{array}{l}\langle \theta, \wt()1( \CY_k(v,x)u)\rangle
=\langle \theta, \wt()1(\sum_m v_{m,k}uz^{-m-1})\rangle \cr
\mbox{}\qquad=(\gr(u)+\gr(v))\langle \theta,\sum_m v_{m,k}uz^{-m-1}\rangle
+\langle \theta, \sum_m (-m-1) v_{m,k}uz^{-m-1}\rangle \cr
\mbox{}\qquad=(\gr(v)+\gr(u))\langle \theta, \CY_k(v,z)u\rangle
+\frac{\pd}{\pd z}\langle \theta, \CY_k(v,z)u\rangle z.
\end{array}$$ 
Using the notation $\Xi=L(0)\!-\!\wt()1$ and 
$(X+Y+Z)^{\ell}=\sum_{i,j,p\in \bN} \lambda^{\ell}_{i,j,p}X^iY^jZ^p$, we obtain  
$$\begin{array}{l}
\langle \theta, (k+1)\CY_{k+1}(v,z)u\rangle \cr
\mbox{}\qquad=\langle \theta, (L(0)\!-\!\wt()1)\CY_k(v,x)u+\CY_k(v,x)(L(0)\!-\!\wt()1)u
+\CY_k((L(0)\!-\!\wt()1)v,x)u\rangle \cr
\mbox{}\qquad=\langle \theta, \Xi\CY_k(v,x)u\rangle
+\langle \theta, \CY_k(v,x) \Xi u\rangle
+\langle \theta, \CY_k(\Xi v,x)u\rangle \cr
\mbox{}\qquad=\cdots=\sum_{\{i,j,p\in \bN\mid i+j+p=\ell\leq k\}} \lambda^{\ell}_{i,j,p}\langle \theta, \Xi^i(\CY_{k+1-\ell}(\Xi^jv,z)\Xi^p u\rangle.
\end{array} $$
Since there is $\bar{K}\in \bN$ such that 
$(\Xi)^{\bar{K}}A=(\Xi)^{\bar{K}}B=(\Xi)^{\bar{K}}C=0$ by Lemma \ref{Nil},  
we have $\langle \theta, \CY_{k}(v,z)u\rangle=0$ 
for $k\geq 3\bar{K}$, that is, we can take 
$K(\theta,v,u)\leq 3\bar{K}$ for any $\theta,v,u$. 
\prend

We will use the following theorem after the proof of Theorem 4.1. 

\begin{thm}\label{Logf}
Let $A,B,C\in \CN\!\CC_1$ and $N\in \bN$. 
For a local intertwining operator 
$\CY$ of type $\binom{C}{A\,\, B}$ over $U$ satisfying ${\rm (I\,\, 1)}\sim {\rm (I\,\, 3)}$, 
there is a finite subset $\Delta_N$ of $\bC$ and $K\in \bN$ such that 
$\langle \theta, \CY(v,z)u\rangle$ has the following expression: 
$$\langle \theta, \CY(v,z)u\rangle=\sum_{t=0}^K\sum_{-s-1\in \Delta_N-\gr(v)-\gr(u)+\gr(\theta)} \langle \theta, v_{s,t}u\rangle z^{-s-1}\log^t(z) $$
for any $\theta \in C^{\vee}_{(\leq N)}, v\in A, u\in B$, where $v_{s,t}\in \Hom(B,C)$ for $v\in V$ and 
we take a principal branch of $z^{-s}$ and $\log(z)$. Furthermore, 
$\widetilde{\CY}(v,z)=\sum_{t=0}^K\sum_{s\in \bC} v_{s,t}z^{-1}\log^t(z)\in I\binom{C}{A\,\,B}$. 
Namely, $\CY$ is a branch of a logarithmic intertwining operator $\widetilde{\CY}$ on $U$.  
\end{thm}

\pr 
Set $S(\theta,v,u;z)=\langle \theta,\CY(v,z)u\rangle$. 
Choose $N\in \bN$ so that $\gr(\theta)\leq N$. 
Let $Q_N=\{\theta^p: p\in \CQ_{N}\}$, $J_A=\{v^i\mid i\in \CP_A\}$, $J_B=\{u^j\mid j\in \CP_B\}$ be bases of $C^{\vee}_{(\leq N)}$, $P_A$, $P_B$, respectively.
From ${\rm (I\,\, 1)}$ and ${\rm (I\,\, 2)}$, for $\alpha\in V_{\geq 1}$, we have 
$$\begin{array}{l}
S(\theta, \alpha_{-1}v,u;z)
=\sum_{j=0}^{\infty} S((\alpha_{-1-j})^{\ast}\theta, v,u;z)\rangle z^j
+\sum_{j=0}^{\infty} S(\theta, v,\alpha_ju;z)z^{-1-j} \quad\mbox{  and}\cr
S(\theta, v,\alpha_{-1}u;z)
=S((\alpha_{-1})^{\ast}\theta, v,u;z)
-\sum_{j=0}^{\infty}(-1)^j S(\theta, \alpha_j v,u;z)z^{-1-j}.
\end{array}$$
Define $\gr^-(\theta,v,u)=\gr(v)+\gr(u)-\gr(\theta)$ and the total grade 
$\gr^+(\theta,v,u)=\gr(v)+\gr(u)+\gr(\theta)$ and 
set $T(\theta,v,u;z)=F(\theta,v,u;z) z^{\gr^-(\theta,v,u)}$. Then we obtain  
$$\begin{array}{l}
T(\theta, \alpha_{-1}v,u;z)=\sum_{j=0}^{\infty} T((\alpha_{-1-j})^{\ast}\theta,v,u;z)
+\sum_{j=0}^{\infty} T(\theta,v,\alpha_ju;z) \quad\mbox{  and} \cr
T(\theta,v,\alpha_{-1}u;z)=T((\alpha_{-1})^{\ast}\theta,v,u;z)
-\sum_{j=0}^{\infty}(-1)^jT(\theta,\alpha_jv,u;z). 
\end{array}\eqno{(2.6)}$$

\begin{cmt}\label{Redu} (1) The coefficients are all constants. \\
(2) The total grades of triples in RHS of (2.6) are less than 
the total grade $\gr(\theta)+\gr(v)+\gr(u)+\gr(\alpha)$ of LHS of (2.6). 
Hence for $(\theta,v,u)$ with $\gr(\theta)\leq N$,  
by iterating the expansion (2.6) as long as we have elements of $C_1(A)$ or $C_1(B)$ 
in the second and third coordinates, 
we finally obtain that $T(\theta,v,u;z)$ is a linear combination of 
$\{T(\theta^p,v^i,u^j;z)\mid (p,i,j)\in \CQ_N\times \CP_A\times \CP_B\}$ with constant 
coefficients.  
\end{cmt}
By the direct calculation, we also have  
$$\begin{array}{l}
\frac{d}{dz}T(\theta,v,u;z)=\frac{d}{dz}\{S(\theta,v,u;z)z^{\gr^-(\theta,v,u)}\}\cr
\mbox{}\qquad=S(\theta, L(-1)v,u;z)z^{\gr^-(\theta,v,u)}
+(\gr^-(\theta,v,u))S(\theta,v,u;z)z^{\gr^-(\theta,v,u)-1}\cr
\mbox{}\qquad=T(\theta, L(-1)v,u;z)z^{-1}+(\gr^-(\theta,v,u))G(\theta,v,u;z)z^{-1}. 
\end{array}$$
Therefore, for a vector valued function 
$Z=(T(\theta^p,v^i,u^j;z))_{(p,i,j)\in \CQ_N\times \CP_A\times \CP_B}$ of size $s$, there is a matrix $A_0\in M_{s\times s}(\bC)$ such that 
$$\frac{d}{dz}Z=\frac{A_0}{z}Z,  \eqno{(2.7)}$$
where $s=|\CQ_N|\times |\CP_A|\times |\CP_B|$, that is, $Z$ satisfies an Euler system (2.7) and 
$T(\theta,v,u;z)$ is a component of a solution of (2.7) for $(\theta,v,u)\in Q_N\times J_A\times J_B$. By Lemma \ref{Jordan}, 
all components of solutions $T(\theta,v,u;z)$ of (2.7) on $\bC\backslash\{0\}$ 
are linear combinations of $\{z^s\log^t(z)\mid s\in \bC, t\in \bN \}$. 
Namely, for each $N,v,u$ with $N\geq \gr(\theta)$, 
there is a finite set $\Delta_N\subseteq \bC$ and $K(N,v,u)\in \bN$ such that  
$$\langle \theta, \CY(v,z)u\rangle=\sum_{s_j\in \Delta_N-\gr(v)-\gr(u)+\gr(\theta)}\sum_{t=0}^{K(N,v,u)}
\langle \theta, v_{-s_j-1,t}u\rangle z^{s_j}\log^t(z) \eqno{(2.8)}$$
with $v_{-s_j-1,t}\in \Hom(B,C)$. 
Since $A$ and $B$ are generated from $P_A$ and $P_B$ as $V$-modules, respectively, we know that 
Eq.(2.8) holds for any $\theta\in C^{\vee}, v\in A, u\in B$ and $z\in U$ 
by using (2.6). 
Although our argument depends on the choice of $N\geq \gr(\theta)$, we note that 
there is $K$ such that 
$K(N,v,u)\leq K$ by Lemma 2.4. 
Therefore it is easy to check that 
$\widetilde{\CY}(v,z)=\sum_{t=0}^K\sum_{s\in \bC} v_{s,t}z^{-s-1}\log^t(z)$ is a 
logarithmic intertwining operator of type $\binom{C}{A\,\, B}$ and 
$\CY$ is a branch of it on $U$. 
\prend

\noindent
{\bf $[$Fusion product$]$} \quad 
In this paper, for $\bN$-gradable $V$-modules $A$ and $B$, a fusion product 
$(A\boxtimes B, \CY^{AB})$ 
is defined to be a pair of an $\bN$-gradable $V$-module $A\boxtimes B$ and $\CY^{AB}\in I\binom{A\boxtimes B}{A\,\, B}$ satisfying the universal property for $\bN$-gradable modules, that is, 
for any $\bN$-gradable $V$-module $C$ and $\CY\in I\binom{C}{A\,\, B}$, there is the unique isomorphism $\phi:A\boxtimes B \to C$ such that 
$\phi(\CY^{AB}(v,z)u)=\CY(v,z)u$ for $v\in A$ and $u\in B$.

\subsection{Borcherds Identities for four-point correlation functions}
From the commutativity ${\rm (I\,\,1)}$ and associativity ${\rm (I\,\,2)}$, we have the following identities: 

\begin{lmm}[Borcherds identities]\label{Cal2}
Let $\CY^1\in I\binom{A(BC)}{A\, \, B\boxtimes C}$, $\CY^2\in I\binom{B\boxtimes C}{B\,\, C}$,
$\CY^{3}\in I\binom{A\boxtimes B}{A\,\,B}$, $\CY^{4}\in I\binom{(AB)C}{A\boxtimes B\,\, C}$. 
For $\theta\in (A(BC))^{\vee}, \theta'\in ((AB)C)^{\vee}, v\in A, u\in B$, and $w\in C$ we set 
$F^{12}(\vec{\xi};x,y)\!=\!\langle \theta, \CY^1(v,x)\CY^2(u,y)w\rangle$ 
and $F^{34}(\vec{\xi}';x,y)\!=\!\langle \theta', \CY^{3}(\CY^{4}(v,x\!-y)u,y)w\rangle$. 
Then for $\alpha\in V$ and $n\in \bZ$, 
as formal $\bC$-power series with logarithm functions, we have:
$$\begin{array}{l}
(1A)\,\, F^{12}({(\alpha_n)^{\ast}}^{[1]}\vec{\xi};x,y)
=\sum_{j=0}^{\infty}\binom{n}{j}F^{12}(\alpha^{[2]}_j\vec{\xi};x,y){x^{n-j}}\cr
\mbox{}\quad \quad +\sum_{j=0}^{\infty}\binom{n}{j}F^{12}(\alpha^{[3]}_j\vec{\xi};x,y)y^{n-j} +F^{12}(\alpha_{n}^{[4]}\vec{\xi};x,y).\cr
(1B)\,\, F^{34}({(\alpha_{n})^{\ast}}^{[1]}\vec{\xi}';x,y)
=\sum_{j=0}^{\infty}\binom{n}{j}
F^{34}(\alpha_j^{[2]}\vec{\xi}';x,y) \iota_{y,x-y}\{x^{n-j}\} \cr
\mbox{}\quad\quad +\sum_{j=0}^{\infty} 
\binom{n}{j}F^{34}(\alpha_j^{[3]}\vec{\xi}';x,y)y^{n-j}+F^{34}(\alpha_{n}^{[4]}\vec{\xi}';x,y).\hspace{57mm}
\end{array}$$ 
$$\begin{array}{l}
(2A) \,\,F^{12}(\alpha_n^{[2]}\vec{\xi};x,y)
=\sum_{j=0}^{\infty}\binom{n}{j}F^{12}({(\alpha_{n-j})^{\ast}}^{[1]}\vec{\xi};x,y){(-x)^j} \cr
\mbox{}\quad -\sum_{j=0}^{\infty} \binom{n}{i}F^{12}(\alpha_i^{[3]}\vec{\xi};x,y)\iota_{x,y}\{(-x+y)^{n-j}\}-\sum_{j=0}^{\infty}\binom{n}{j}F^{12}(\alpha_j^{[4]}\vec{\xi};x,y){(-x)^{n-j}}.\cr 
(2B) \,\,F^{34}(\alpha_n^{[2]}\vec{\xi}';x,y)=\!\sum_{j=0}^{\infty}
\binom{n}{j}F^{34}({(\alpha_{n-j})^{\ast}}^{[1]}\vec{\xi}';x,y) 
 \iota_{y,x-y}\{(-x)^{j}\}\cr
\mbox{}\quad -\sum_{j=0}^{\infty}\binom{n}{j}F^{34}(\alpha_j^{[3]}\vec{\xi}';x,y){(-x\!+\!y)^{n-j}}-\sum_{j=0}^{\infty}\binom{n}{j}F^{34}(\alpha_j^{[4]}\vec{\xi}';x,y)\iota_{y,x-y}\{(-x)^{n-j}\}.\hspace{4mm}
\end{array}$$
$$\begin{array}{l}
(3A)\,\, F^{12}(\alpha_{n}^{[3]}\vec{\xi};x,y) 
=\sum_{j=0}^{\infty}\binom{n}{j}(-1)^jF^{12}({(\alpha_{n-j})^{\ast}}^{[1]}\vec{\xi};x,y)
{y^j} \cr
\mbox{}\quad \quad-\sum_{j=0}^{\infty}\binom{n}{j}F^{12}(\alpha_j^{[2]}\vec{\xi};x,y) \iota_{x,y}\{ (x-y)^{n-j}\}-\sum_{j=0}^{\infty}\binom{n}{j}(-1)^{j+n}F^{12}(\alpha_i^{[4]}\vec{\xi};x,y) {y^{n-j}}.  \cr
(3B)\,\, F^{34}(\alpha_{n}^{[3]}\vec{\xi}';x,y) 
\!=\!\sum_{j=0}^{\infty}\binom{n}{j}(-1)^jF^{34}({(\alpha_{n-j})^{\ast}}^{[1]}\vec{\xi}';x,y)
{y^j} \cr
\mbox{}\quad \quad-\sum_{j=0}^{\infty}\binom{n}{j}F^{34}(\alpha_j^{[2]}\vec{\xi}';x,y)
{(x-y)^{n-j}}-\sum_{j=0}^{\infty}\binom{n}{j}(-1)^{j+n}
F^{34}(\alpha_j^{[4]}\vec{\xi}';x,y){y^{n-j}}.
\end{array}$$
$$\begin{array}{l}
(4A)\,\, F^{12}(\alpha_{n}^{[4]}\vec{\xi};x,y) 
=F^{12}({(\alpha_{n})^{\ast}}^{[1]}\vec{\xi};x,y)-\sum_{j=0}^{\infty} \binom{n}{j}F^{12}(\alpha_j^{[2]}\vec{\xi};x,y)
{x^{n-j}}\cr
\mbox{}\quad \quad -\sum_{j=0}^{\infty}\binom{n}{j}F^{12}(\alpha_j^{[3]}\vec{\xi};x,y){y^{n-j}}. \cr
(4B)\,\, F^{34}(\alpha_n^{[4]}\vec{\xi}';x,y) 
=F^{34}({(\alpha_{n})^{\ast}}^{[1]}\vec{\xi}';x,y)-\sum_{j=0}^{\infty}\binom{n}{j}
F^{34}(\alpha_j^{[2]}\vec{\xi}';x,y)\iota_{y,x-y}\{x^{n-j}\}\hspace{6mm}\cr
\mbox{}\quad \quad -\sum_{j=0}^{\infty} \binom{n}{j}F^{34}(\alpha_j^{[3]}\vec{\xi}';x,y)
{y^{n-j}}. 
\end{array}$$
Here $(\alpha_m)^{\ast}$ denotes its adjoint operator of $\alpha_m$ and 
$\iota_{x,y}\{(x-y)^s\}:=\sum_{j=0}^{\infty}\binom{s}{j}x^{s-j}(-y)^j$ and 
$\iota_{y,x-y}\{x^s\}:=\iota_{y,x-y}\{(y+(x-y))^s\}=\sum_{j=0}^{\infty}\binom{s}{j}y^{s-j}(x-y)^j$.
\end{lmm}

\begin{cmt}\label{Same} 
(1) The coefficients are in $\bC[x^{\pm 1},y^{\pm 1}, \tau_{x,y}\{(x-y)^{\pm 1}\}, 
\tau_{y,x-y}\{x^{\pm 1}\}]$. Hence 
the same equations also hold for coefficients of logarithm functions. \\
(2) Then for $n\in \bZ$, $m,\ell=1,2,3,4$, and $j\geq 0$,  
the coefficients of $F^{12}(\alpha^{[\ell]}_j\vec{\xi};x,y)$ 
in RHS of $(mA)$ for 
$F^{12}(\alpha^{[m]}_n\vec{\xi};x,y)$ are equal to 
the corresponding coefficients of $F^{34}(\alpha^{[\ell]}_j\vec{\xi};x-y,y)$ in RHS of 
$(mB)$ for $F^{34}(\alpha^{[m]}_n\vec{\xi};x-y,y)$ 
excepting the expansions by $\iota_{x,y}$ or $\iota_{x-y,y}$.
Hence the corresponding coefficients in (mA) and (mB) are same on a domain 
$\{(x,y)\in \bC^2 \mid 0<|x-y|<|y|<|x|\}$. 
\end{cmt} 

\pr 
First, we prepare the following equations: 
$$\begin{array}{lr}
\sum_{\ell=0}^{\infty}\sum_{j=0}^{\infty}\binom{n}{j}
\binom{j}{\ell}x^{n-j}y^{j-\ell}Z^{\ell}=
\sum_{\ell=0}^{\infty}\binom{n}{i}\iota_{x,y}\{(x+y)^{n-\ell}\}Z^{\ell} 
&(2.9)\cr
\sum_{i,j\in \bN}\binom{n-j}{i}\binom{n}{j}(-1)^{i+j}(x-y)^jy^iZ^{i+j}
=\sum_{\ell=0}^{\infty}\binom{n}{\ell}(-1)^{\ell}\iota_{y,x-y}\{x^{\ell}\}Z^{\ell} 
\hspace{2cm}&(2.10)
\end{array}$$
and we will use the equations between the coefficients of $Z^{\ell}$ of LHS and RHS. 

By commutativity ${\rm (I\,\,1)}$ and associativity ${\rm (I\,\,2)}$, we have:

\noindent
$$\begin{array}{rl}
(1A)&\langle(\alpha_n)^{\ast}\theta, \CY^1(v,x)\CY^2(u,y)w\rangle 
=\langle \theta, \alpha_n\CY^1(v,x)\CY^2(u,y)w\rangle 
\hspace{2.5cm} \hfill \cr 
=&\langle \theta, \CY^1(v,x)\CY^2(u,y)(\alpha_{n}w)\rangle+\sum_{j=0}^{\infty}\binom{n}{j}\langle \theta, \CY^1(\alpha_jv,x)\CY^2(u,y)w\rangle{x^{n-j}}\cr
&+\sum_{j=0}^{\infty}\binom{n}{j}\langle \theta, \CY^1(v,x)\CY^2(\alpha_ju,y)w\rangle 
{y^{n-j}}. \cr
(1B) &\langle \alpha_{a}\theta', \CY^{3}(\CY^{4}(v,x-y)u,y)w\rangle =
\langle \theta', \alpha_n\CY^{3}(\CY^{4}(v,x-y)u,y)w\rangle \hspace{1cm} \hfill\cr
=&\langle \theta', \CY^{3}(\CY^{4}(v,x-y)u)y)\alpha_{n}w\rangle+\sum_{j=0}^{\infty} \binom{n}{j}\langle \theta',\CY^3 ( \alpha_j\CY^4(v,x-y)u),y)w\rangle \cr
=&\langle \theta', \CY^{3}(\CY^{4}(v,x-y)u,y)\alpha_nw\rangle +\sum_{j=0}^{\infty} \binom{n}{j}\langle \theta', \CY^{3}(\CY^{4}(v,x-y)\alpha_ju,y)w\rangle 
{y^{n-j}}\cr
&+\sum_{j=0}^{\infty}\binom{n}{j}\sum_{i=0}^{\infty}\binom{j}{i}
\langle \theta', \CY^{3}(\CY^{4}((\alpha_iv),x-y)u),y)w\rangle {(x-y)^{j-i}y^{n-j}} \cr
=&\langle \theta', \CY^{3}(\CY^{4}(v,x-y)u)y)\alpha_{n}w\rangle
+\sum_{j=0}^{\infty} \binom{n}{j}\langle \theta', \CY^{3}(\CY^{4}(v,x-y)\alpha_ju,y)w\rangle {y^{n-j}}\cr
&+\sum_{j=0}^{\infty}\binom{n}{j}
\langle \theta', \CY^{3}(\CY^{4}(\alpha_iv,x-y)u,y)w\rangle {\iota_{y,x-y}\{x^{n-j}\}}  
\mbox{  by $(2.9)$}. 
\end{array}$$
$$\begin{array}{rl}
(2A)& \langle \theta, \CY^1(\alpha_nv, x)\CY^2(u,y)w\rangle=\sum_{j=0}^{\infty}\binom{n}{j}(-1)^j\langle \theta, \alpha_{n-j}
\CY^1(v,x)\CY(u,y)w\rangle {x^j}\cr
&-\sum_{j=0}^{\infty}\binom{n}{j}(-1)^{j+n}\langle \theta,\CY^1(v,x)\alpha_j\CY^2(u,y)w\rangle x^{n-j}\cr
=&\sum_{j=0}^{\infty}\binom{n}{j}(-1)^j\langle (\alpha_{n-j})^{\ast}\theta, \CY^1(v,x)\CY^2(u,y)w\rangle  
{x^j} \cr
&-\sum_{j=0}^{\infty}\binom{n}{j}(-1)^{j+n}\langle \theta,\CY^1(v,x)\CY^2(u,y)(\alpha_jw)\rangle 
{x^{n-j}}\cr
&-\sum_{j=0}^{\infty}\binom{n}{j}(-1)^{j+n}\sum_{i=0}^{\infty}\binom{j}{i}\langle \theta,\CY^1(v,x)\CY^2(\alpha_iu,w)
\rangle {x^{n-j}y^{j-i}} \cr
=&\sum_{j=0}^{\infty}\binom{n}{j}(-1)^{j}\langle (\alpha_{n-j})^{\ast}\theta, \CY^1(v,x)\CY^2(u,y)w\rangle  
{x^j} \cr
&-\sum_{j=0}^{\infty}\binom{n}{j}(-1)^{n+j}\langle \theta,\CY^1(v,x)\CY^2(u,y)(\alpha_jw)\rangle 
{x^{n-j}}\cr
&-\sum_{j=0}^{\infty} \binom{n}{j}(-1)^{n+j}\langle \theta,\CY^1(v,x)\CY^2(\alpha_ju,w)
\rangle {\iota_{x,y}\{(x-y)^{n-j}\}} \mbox{  by $(2.9)$}. \cr
(2B)& \langle \theta', \CY^3(\CY^4(\alpha_nv,x-y)u,y)w\rangle\cr
=&\sum_{j=0}^{\infty}\binom{n}{j}(-1)^j \langle \theta', 
\CY^3(\alpha_{n-j}\CY^4(v,x-y)u,y)w\rangle (x-y)^j\cr
&-\sum_{j=0}^{\infty}\binom{n}{j}(-1)^{j+n}\langle \theta', \CY^3(\CY^4(v,x-y)\alpha_ju,y)w\rangle (x-y)^{n-j} \cr
=&\sum_{j=0}^{\infty}\sum_{i=0}^{\infty}\binom{n-j}{i}\binom{n}{j}(-1)^{i+j}
\langle\theta',\alpha_{n-i-j}\CY^3(\CY^4(v,x-y)u,y)w\rangle 
 {(x-y)^jy^i}\cr
&-\sum_{j=0}^{\infty}\sum_{i=0}^{\infty}\binom{n}{j}\binom{n-j}{i}(-1)^{i+n}\langle \theta', \CY_h^3(\CY_k^4(v,x-y)u,y)\alpha_iw\rangle   {(x-y)^jy^{n-i-j}}\cr
&-\sum_{j=0}^{\infty}\binom{n}{j}(-1)^{j+n}\langle \theta', \CY^3(\CY^4(v,x-y)\alpha_ju,y)w\rangle 
 {(x-y)^{n-j}}\cr
=&\sum_{j=0}^{\infty}\binom{n}{j}(-1)^{j}
\langle(\alpha_{n-j})^{\ast}\theta',\CY^3(\CY^4(v,x-y)u,y)w\rangle 
 {\iota_{y,x-y}\{x^j\}}  \mbox{  by $(2.10)$}\cr
&-\sum_{j=0}^{\infty}\binom{n}{j}(-1)^{n+j}\langle \theta', \CY^3(\CY^4(v,x-y)u,y)\alpha_jw\rangle   {\iota_{y,x-y}\{x^{n-j}\}} \mbox{  by $(2.9)$}\cr
&-\sum_{j=0}^{\infty}\binom{n}{j}(-1)^{j+n}\langle \theta', \CY^3(\CY^4(v,x-y)\alpha_ju,y)w\rangle 
 {(x-y)^{n-j}}. 
\end{array}$$
$$\begin{array}{rl}
(3A)&\langle \theta, \CY^1(v,x)\CY^2(\alpha_{n}u,y)w\rangle
=\sum_{i=0}^{\infty}\binom{n}{i}(-1)^i \langle \theta, \CY^1(v,x)\alpha_{n-i}\CY^2(u,y)w\rangle y^j\cr
&+\sum_{i=0}^{\infty} \binom{n}{i}(-1)^{i+n+1} 
\langle \theta, \CY^1(v,x)\CY^2(u,y)\alpha_iw\rangle y^{n-i} \cr
=&\sum_{i=0}^{\infty}\binom{n}{i}(-1)^i\langle \theta, \alpha_{n-i}\CY^1(v,x)\CY^2(u,y)w\rangle {y^i}\cr
&-\sum_{i=0}^{\infty}\sum_{j=0}^{\infty}\binom{n}{i}(-1)^i\binom{n-i}{j}\langle \theta, \CY^1(\alpha_jv,x)\CY^2(u,y)w\rangle {x^{n-i-j}y^i} \cr
&-\sum_{i=0}^{\infty}\binom{n}{i}(-1)^{i+n}\langle \theta, \CY^1(v,x)\CY^2(u,y)(\alpha_iw)\rangle {y^{n-i}} \cr
=&\sum_{j=0}^{\infty}\binom{n}{j}(-1)^j\langle (\alpha_{n-j})^{\ast}\theta, \CY^1(v,x)\CY^2(u,y)w\rangle {y^j}\cr
&-\sum_{j=0}^{\infty}\binom{n}{j}\langle \theta, \CY^1(\alpha_jv,x)\CY^2(u,y)w\rangle {\iota_{x,y}\{(x-y)^{n-j}\}} \mbox{  by $(2.9)$}\cr
&-\sum_{j=0}^{\infty}\binom{n}{j}(-1)^{j+n}\langle \theta, \CY^1(v,x)\CY^2(u,y)\alpha_iw\rangle {y^{n-j}}. \cr
(3B)&\langle \theta', \CY^{3}(\CY^{4}(v,x\!-\!y)(\alpha_{n}u),y)w\rangle
=\langle \theta', \CY^3(\alpha_{n}\CY^4(v,x-y)u,y)w\rangle \cr
 &-\sum_{j=0}^{\infty}\binom{n}{j}\langle \theta, \CY^3(\CY^4(\alpha_jv,x-y)u,y)w
\rangle (x-y)^{n-j} \cr
=&\sum_{j=0}^{\infty}\binom{n}{j}(-1)^j\langle 
(\theta', \alpha_{n-j}\CY^3(\CY^4(v,x-y)u,y)w\rangle y^j\cr
&-\sum_{j=0}^{\infty}\binom{n}{j}(-1)^{j+n}\langle \theta', \CY^3(\CY^4(v,x-y)u,y)\alpha_jw)\rangle y^{n-j} \cr
&-\sum_{j=0}^{\infty}\binom{n}{j}\langle \theta', \CY^3(\CY^4(\alpha_jv,u,y)w\rangle
(x-y)^{n-j}\cr
=&\sum_{j=0}^{\infty}\binom{n}{j}(-1)^j\langle (\alpha_{n-j})^{\ast}\theta', \CY^{3}(\CY^{4}(v,x\!-\!y)u,y)w\rangle {y^j}\cr
&-\sum_{j=0}^{\infty}\binom{n}{j}(-1)^{j+n}\langle \theta', \CY^{3}(\CY^{4}(v,x\!-\!y)u,y)\alpha_jw\rangle {y^{n-j}}\cr
&-\sum_{j=0}^{\infty}\binom{n}{j}\langle \theta', \CY^{3}(\CY^{4}(\alpha_jv,x\!-\!y)u,y)w\rangle {(x-y)^{n-j}}.
\end{array}$$
$$\begin{array}{rl}
(4A)&\langle \theta, \CY^1(v,x)\CY^2(u,y)(\alpha_nw)\rangle
=\langle (\alpha_n)^{\ast}\theta, \CY^1(v,x)\CY^2(u,y)w\rangle \cr
&-\sum_{j=0}^{\infty}\binom{n}{j}\langle \theta, \CY^1(v,x)\CY^2(\alpha_ju,y)w\rangle {y^{n-j}}
\!-\!\sum_{j=0}^{\infty} \binom{n}{j}\langle \theta, \CY^1(\alpha_jv,x)\CY^2(u,y)w\rangle {x^{n-j}}.\cr
(4B)&\langle \theta', \CY^{3}(\CY^{4}(v,x-y)u,y)(\alpha_nw)\rangle\cr
=&\langle \theta', \alpha_n\CY^3(\CY^4(v,x-y)u,y)w\rangle-\sum_{j=0}^{\infty}\binom{n}{j}\langle \theta', \CY^3(\alpha_j(\CY^4(v,x-y)u,y)w\rangle\cr
=&\langle (\alpha_n)^{\ast}\theta', \CY^3(\CY^4(v,x-y)u),y)w\rangle 
\!-\!\sum_{j=0}^{\infty} \binom{n}{j}\langle \theta', \CY^3(\CY^4(v,x-y)(\alpha_ju)),y)w\rangle {y^{n-j}}\cr
&\!\!-\!\sum_{j=0}^{\infty}\sum_{i=0}^{\infty}\binom{n}{j}\binom{j}{i}\langle \theta', \CY^{3}(\CY^4(\alpha_iv,x-y)u,y)w\rangle {(x-y)^{j-i}y^{n-j}}\cr
=&\langle (\alpha_n)^{\ast}\theta', \CY^{3}(\CY^{4}(v,x-y)u),y)w\rangle 
-\sum_{j=0}^{\infty} \binom{n}{j}\langle \theta', \CY^{3}(\CY^{4}(v,x-y)\alpha_ju,y)w\rangle {y^{n-j}}\cr
&\!\!-\!\sum_{i=0}^{\infty}\binom{n}{i}\langle \theta', \CY^{3}(\CY^4(\alpha_iv,x-y)u,y)w\rangle {\iota_{y,x-y}\{x^{n-i}\}} \mbox{  by $(2.9)$}. \qed 
\end{array}$$

Let us explain a role of $L(-1)$-derivative properties ${\rm (I\,\, 3)}$. 

\begin{rmk}\label{Derivation}
From $L(-1)$-derivative property given by $L(-1)^{[3]}$, 
we have: \vspace{-2mm}
$$\begin{array}{rl}
\frac{\pd}{\pd y}(\sum F^{A(BC)}_{h,k}(\vec{\xi};x,y)\log^h(x)\log^k(y))&=\frac{\pd}{\pd y} F^{A(BC)}(\vec{\xi};x,y)=
F^{A(BC)}(L(-1)^{[3]}\vec{\xi};x,y)\cr
=&\sum_{h,k}F_{h,k}^{A(BC)}(L(-1)^{[3]}\vec{\xi};x,y)\log^h(x)\log^k(y).
\end{array}\vspace{-2mm}$$ 
By comparing the coefficients of $\log^h(x)\log^k(y)$, we get \vspace{-3mm} 
$$F^{A(BC)}_{h,k}(L(-1)^{[3]}\vec{\xi};x,y)=\frac{\pd}{\pd y}
F^{A(BC)}_{h,k}(\vec{\xi};x,y)+\frac{k+1}{y} F^{A(BC)}_{h,k+1}(\vec{\xi};x,y). 
\eqno{(2.11)}\vspace{-3mm}$$
Similarly, by considering $\frac{\pd}{\pd x}$ and $L(-1)^{[2]}$, we obtain: \vspace{-3mm}
$$F^{A(BC)}_{h+1,k}(\vec{\xi};x,y)=\frac{x}{h+1} \{ F^{A(BC)}_{h,k}(L(-1)^{[2]}\vec{\xi};x,y)
-\frac{\pd}{\pd x} F^{A(BC)}_{h,k}(\vec{\xi};x,y) \}. \eqno{(2.12)}\vspace{-3mm}$$ 
In particular, we can recover all $F^{A(BC)}_{h,k}$ 
from $F^{A(BC)}_{0,0}$ by $\frac{\pd}{\pd x}$ and $\frac{\pd}{\pd y}$.  
\end{rmk}

\section{The proof of theorems 3.1, 3.2 and Corollary 3.3}
\subsection{Theorem 3.1 (Reduction to small weight subspaces)}
We recall Lemma 2.11. See Table 1 for Borcherds identities $(2A)\sim(4A)$ for 
$n=-1$.   
\begin{table}[hbtp]
  \caption{Borcherds identity}
  \label{table:data_type}
  \begin{tabular}{|ll|l|l|}
\hline
  Borcherds identities &D:coefficient& q &$D\times (x-y)^q\mbox{ and }D\times y^q$ \\
\hline
$(2A)$ $F(\alpha_{-1}^{[2]}\vec{\xi})\times$ &1 & 0 & 1\\
$\mbox{}\quad=\sum_{j=0}^{\infty}F(((\alpha_{-1-j})^{\ast})^{[1]}\vec{\xi})\times$ &$x^j$ & $\wt(\alpha)$ &no negative powers \\
$\mbox{}\qquad+\sum_{j=0}^{\infty}F(\alpha_j^{[3]}\vec{\xi})\times$ &
$\hspace{-6mm}\iota_{x,y}\{(x-y)^{-j-1}\}$& j+1 &$1$ and $(x-y)^{-j-1}y^{j+1}$\\
$\mbox{}\qquad+\sum_{j=0}^{\infty}F(\alpha_j^{[4]}\vec{\xi})\times$ &${x^{-j-1}}.$ &j+1&$x^{-j-1}(x\!-\!y)^{j+1}$ and $(y/x)^{j+1}$\\
$(3A)$ $F(\alpha_{-1}^{[3]}\vec{\xi})\times$ &1 & 0 & 1\\
$\mbox{}\quad=\sum_{j=0}^{\infty}F(((\alpha_{-1-j})^{\ast})^{[1]}\vec{\xi})\times$ &$ {y^j}$ &
$\wt(\alpha)$&no negative powers\\
$\mbox{}\qquad-\sum_{j=0}^{\infty}(-1)^jF( \alpha_j^{[2]}\vec{\xi})\times$ 
&$\hspace{-6mm}{\iota_{x,y}\{(x-y)^{-j-1}\}}$ &$j+1$&$1$ and $(x-y)^{-j-1}y^{j+1}$\\
$\mbox{}\qquad+\sum_{j=0}^{\infty}F(\alpha_j^{[4]}\vec{\xi})\times$ &$ {y^{-j-1}}.$ &$ j+1$ &$y^{-j-1}(x-y)^{j+1}$ and $1$\\
$(4A)$ $F(\alpha_{-1}^{[4]}\vec{\xi})\times$ &1 & 0 & 1\\
$\mbox{}\quad=F(((\alpha_{-1})^{\ast})^{[1]}\vec{\xi})\times$ &1 & $\wt(\alpha)$ 
&no negative powers\\
$\mbox{}\qquad-\sum_{j=0}^{\infty}(-1)^j F(\alpha_j^{[3]}\vec{\xi})\times$ 
& $y^{-j-1}$ &j+1&$(x-y)^{j+1}y^{-j-1}$ and $1$\\
$\mbox{}\qquad-\sum_{j=0}^{\infty}(-1)^jF(\alpha_j^{[2]}\vec{\xi})\times$ &${x^{-j-1}}. $&j+1&$(x-y)^{j+1}x^{-j-1}$ and $(y/x)^{j+1}$\\
\hline 
  \end{tabular}\\
\begin{tabular}{|l|}
Here the first columns are Borcherds identities for $n=-1$ and 
$F(\vec{\xi})$ denotes $\mbox{}\hspace{15mm}$\\
$F^{A(BC)}_{h,k}(\vec{\xi};x,y)$ and 
$q=\gr^{234}(\alpha_{-1}^{[m]}\vec{\xi})-\gr^{234}(\alpha_j^{[i]}\vec{\xi})$ for a quadruple $\alpha_j^{[i]}\vec{\xi}$ in (mA). \\
\hline
\end{tabular}\\
\end{table}

\begin{cmt}\label{Redu2} Since $\gr((\alpha_{-1-j})^{\ast}\theta)< \gr(\theta)$, 
we have $(\alpha_{-1-j})^{\ast}\theta\in (A(BC))^{\vee}_{(< N)}$ for $\alpha\in V_{\geq 1}$ and $j\in \bN$. 
Since $\gr(\alpha_jw)=\gr(\alpha_{-1}w)-j-1$ for $j\in \bN$, 
the total grades of quadruples of each term in RHS are less than the total grade of quadruples 
in LHS, that is, 
$\gr(((\alpha_{-1-j})^{\ast})^{[1]}\vec{\xi}), \gr(\alpha_j^{[k]}\vec{\xi}) <\gr(\alpha_{-1}^{[i]}\vec{\xi})$.  
Therefore by iterating these reductions as long as we have an element of $C_1(A), C_1(B)$ or $C_1(C)$ 
in the second, third or fourth coordinates of quadruples, we finally get an expression of 
$F_{h,k}^{A(BC)}(\vec{\xi};x,y)$ as a linear combination of 
$$\{F_{h,k}^{A(BC)}(\vec{\mu};x,y):\vec{\mu}\in \CJ^0_{N} \}$$ with 
coefficients in $\bC[x,x^{-1},y,y^{-1},\iota_{x,y}\{(x-y)^{-1}\}]$
for any $\vec{\xi}\in \Omega_N$. 
\end{cmt}

\begin{cmt}\label{Nega(x-y)}
Since Borcherds identities for $F^{A(BC)}_{h,k}$ and $F^{(AB)C}_{h,k}$ have 
the same coefficients except the way of expansions, 
we can get Borcherds identities $(2B)\sim(4B)$ by viewing 
$F(\vec{\xi})$ as $F^{(AB)C}_{h,k}(\vec{\xi}';x,y)$ and expanding all coefficients 
by $\iota_{y,x-y}$. 
\end{cmt}

\begin{cmt}\label{Nega}
In the case of $D\times (x-y)^q$, negative powers appear only on $x$ or $y$. 
In the case of $D\times y^q$, negative powers occur only on $(x-y)$ or $x$.  
We hence consider 
$$\begin{array}{ll}
G_{h,k}^{A(BC):y}(\vec{\xi};x_0,y)&:=F_{h,k}^{A(BC)}(\vec{\xi};x_0,y)y^{\gr^{234}(\vec{\xi})}, \cr
G^{A(BC):x-y}_{h,k}(\vec{\xi};x,y_0)&:=F^{A(BC)}_{h,k}(\vec{\xi};x,y_0)(x-y_0)^{\gr^{234}(\vec{\xi})}, \cr
G_{h,k}^{(AB)C:y}(\vec{\xi}',w;x_0,y)&:=F_{h,k}^{(AB)C}(\vec{\xi}';x_0,y)y^{\gr^{234}(\vec{\xi})}, \mbox{ and}\cr
G^{(AB)C:x-y}_{h,k}(\vec{\xi}';x,y_0)&:=F^{(AB)C}_{h,k}(\vec{\xi}';x,y_0)(x-y_0)^{\gr^{234}(\vec{\xi})}.
\end{array}
\eqno{(3.1)}$$
\end{cmt}

Then we have:

\begin{thm}\label{G}
For $\vec{\xi}=(\theta,v,u,w)$ with $\theta\in J_N$, $\alpha\in V$, and $x_0\not=0\not=y_0$, \\
(1) $G_{h,k}^{A(BC):y}(\vec{\xi};x_0,y)$ is a linear combination of 
$\{G_{h,k}^{A(BC):y}(\vec{\delta};x_0,y)\mid \vec{\mu}\in 
\CJ^0_N\}$ 
with coefficients 
in $\bC[\iota_{x_0,y}\{(x_0-y)^{-1}\}] [y]\subseteq \bC[[y]]$ and \\
(2) $G_{h,k}^{A(BC):x-y}(\vec{\xi};x,y_0)$ 
is a linear combination of 
$\{G_{h,k}^{A(BC):x-y}(\vec{\mu};x,y_0)\mid 
\vec{\mu}\in \CJ^0_N\}$ with coefficients in $\bC[x,\iota_{y_0,x-y_0}\{x^{-1}\}][x-y_0]\subseteq \bC[[x-y_0]]$. \\
Furthermore, for the residue classes of coefficients modulo $\bC[[x-y_0]](x-y_0)$, we have:
$$\begin{array}{ll}
G_{h,k}^{A(BC):x-y}(\alpha_{-1}^{[2]}\vec{\xi};x,y_0)\equiv \sum_{j=0}^{\infty} G_{h,k}^{A(BC):x-y}(\alpha_j^{[3]}\vec{\xi};x,y_0)
&\pmod{\bC[[x\!-\!y_0]](x\!-\!y_0)},\cr
G_{h,k}^{A(BC):x-y}(\alpha_{-1}^{[3]}\vec{\xi};x,y_0)\equiv \sum_{j=0}^{\infty}(-1)^{j+1} G_{h,k}^{A(BC):x-y}( \alpha_j^{[2]}\vec{\xi};x,y_0) \hspace{-6mm}
&\pmod{\bC[[x\!-\!y_0]](x\!-\!y_0)},\cr
G_{h,k}^{A(BC):x-y}(\alpha_{-1}^{[4]}\vec{\xi};x,y_0)\equiv 0 &\pmod{\bC[[x\!-\!y_0]](x\!-\!y_0)}, \cr
G_{h,k}^{A(BC):y}(\alpha_{-1}^{[2]}\vec{\xi};x_0,y)\equiv 0 &\pmod{\bC[[y]]y},\cr
G_{h,k}^{A(BC):y}(\alpha_{-1}^{[3]}\vec{\xi};x_0,y)\equiv \sum_{j=0}^{\infty} G_{h,k}^{A(BC):y}(\alpha_j^{[4]}\vec{\xi};x_0,y)
&\pmod{\bC[[y]]y},\cr
G_{h,k}^{A(BC):y}(\alpha_{-1}^{[4]}\vec{\xi};x_0,y)\equiv -\sum_{j=0}^{\infty}(-1)^j G_{h,k}^{A(BC):y}( \alpha_j^{[3]}\vec{\xi};x_0,y) 
&\pmod{\bC[[y]]y}.
\end{array}$$
We also obtain the same results for 
$G_{h,k}^{(AB)C:x-y}(\vec{\xi}';x-y_0,y_0)$ and 
$G_{h,k}^{(AB)C:y}(\vec{\xi}';x_0-y,y)$ by 
replacing $A(BC), \vec{\xi}, \vec{\mu}, \CJ_N, \CJ_N^0, \iota_{x,y}$ by 
$(AB)C, \vec{\xi}', \vec{\mu}', \CJ'_N, {\CJ'}_N^0, \iota_{y,x-y}$, respectively. 
\end{thm}

\pr 
We will show the proofs for $G^{A(BC):*}_{h,k}(\vec{\xi};x,y)$. 
To simplify the notation, we denote $G^{A(BC):*}_{h,k}(\vec{\xi};x,y)$ by 
$G^{*}(\vec{\xi};x,y)$ for $\ast=y, x-y$. From Table 1, we have 
$$\begin{array}{l}
G^y(\alpha_{-1}^{[2]}\vec{\xi};x_0,y)=F(\alpha_{-1}^{[2]}\vec{\xi};x_0,y)y^{p+\wt(\alpha)}
=\sum_{j=0}^{\infty} G^y(((\alpha_{-1-j})^{\ast})^{[1]}\vec{\xi};x_0,y)x_0^jy^{\wt(\alpha)}\cr
\mbox{}\quad +\sum_{j=0}^{\infty}G^y(\alpha_j^{[3]}\vec{\xi};x_0,y)x_0\iota_{x_0,y}\{(x_0-y)^{-j-1}\}y^{j+1}
+\sum_{j=0}^{\infty}G^y(\alpha_j^{[4]}\vec{\xi};x_0,y)x_0^{-j-1}y^{j+1}, \cr
G^y(\alpha_{-1}^{[3]}\vec{\xi};x_0,y)=F(\alpha_{-1}^{[3]}\vec{\xi};x_0,y)y^{p+\wt(\alpha)}
=\sum_{j=0}^{\infty} G^y(((\alpha_{-1-j})^{\ast})^{[1]}\vec{\xi};x_0,y) y^{j+\wt(\alpha)}\cr
\mbox{}\quad -\sum_{j=0}^{\infty}(-1)^jG^y(\alpha_j^{[2]}\vec{\xi};x_0,y) \iota_{x_0,y}\{(x_0\!-\!y)^{-j-1}\}y^{j+1}
+\sum_{j=0}^{\infty}G^y(\alpha_j^{[4]}\vec{\xi};x_0,y),  \cr
G^y(\alpha_{-1}^{[4]}\vec{\xi};x_0,y)=F(\alpha_{-1}^{[4]}\vec{\xi};x_0,y)y^{p+\wt(\alpha)}
=G^y(((\alpha_{-1})^{\ast})^{[1]}\vec{\xi};x_0,y)y^{\wt(\alpha)}\cr
\mbox{}\quad -\sum_{j=0}^{\infty}(-1)^jG^y(\alpha_j^{[3]}\vec{\xi};x_0,y)-\sum_{j=0}^{\infty} (-1)^jG^y(\alpha_j^{[2]}\vec{\xi};x_0,y)x_0^{-j-1}y^{j+1}.
\end{array}$$ 
For $y=y_0\not=0$ and a variable $x$, we have:
$$\begin{array}{l}
G^{x-y}(\alpha_{-1}^{[2]}\vec{\xi};x,y_0)=\sum_{j=0}^{\infty} G(((\alpha_{-1-j})^{\ast})^{[1]}\vec{\xi};x,y_0)x^j(x-y_0)^{\wt(\alpha)}\cr
\mbox{}\quad+\sum_{j=0}^{\infty} G(\alpha_j^{[3]}\vec{\xi};x,y_0)\iota_{y_0,x-y_0}\{y_0+(x-y_0)\}\cr
\mbox{}\quad +\sum_{j=0}^{\infty} G(\alpha_j^{[4]}\vec{\xi};x,y_0)\iota_{y_0,x-y_0}\{(y_0+(x-y_0))^{-j-1}\}(x-y_0)^{j+1}, \cr
G^{x-y}(\alpha_{-1}^{[3]}\vec{\xi};x,y_0)=\sum_{j=0}^{\infty} G^{x-y}(((\alpha_{-1-j})^{\ast})^{[1]}\vec{\xi};x,y_0) y_0^{j}(x-y_0)^{\wt(\alpha)}\cr
\mbox{}\quad-\sum_{j=0}^{\infty}(-1)^j G^{x-y}( \alpha_j^{[2]}\vec{\xi};x,y_0) 
+\sum_{j=0}^{\infty}G^{x-y}(\alpha_j^{[4]}\vec{\xi};x,y_0) y_0^{-j-1}(x-y_0)^{j+1}, \cr
G^{x-y}(\alpha_{-1}^{[4]}\vec{\xi};x,y_0)=G^{x-y}(((\alpha_{-1})^{\ast})^{[1]}\vec{\xi};x,y_0)(x-y_0)^{\wt(\alpha)}\cr
\mbox{}\quad-\sum_{j=0}^{\infty}(-1)^jG^{x-y}(\alpha_j^{[3]}\vec{\xi};x,y_0)y_0^{-j-1}(x-y_0)^{j+1}\cr
\mbox{}\quad -\sum_{j=0}^{\infty} (-1)^jG^{x-y}(\alpha_j^{[2]}\vec{\xi};x,y_0)\iota_{y_0,x-y_0}\{(y_0+(x-y_0))^{-j-1}(x-y_0)^{j+1}.
\end{array}$$ 
In particular, for the residue classes of coefficients modulo $\bC[[x-y_0]](x-y_0)$, 
it is easy to check the statements.
Since these processes don't depend on the choice of $G^{A(BC)}$ nor 
$G^{(AB)C}$, we get the same expressions for $G^{(AB)C:y}(\vec{\xi}';x_0,y)$ 
and $G^{(AB)C:x-y}(\vec{\xi}';x,y_0)$
except the notation and expansions by $\iota_{\ast,\ast}$. 

This completes the proof of Theorem \ref{G}.

\subsection{Proof of Theorem 3.2 (Differential systems)}
We will prove the Fuchsian systems (D1) and (D2) for $G^{A(BC):y}_{h,k}$ and 
$G^{A(BC):x-y}_{h,k}$.  
The proofs for the other cases are similar. 
From the $L(-1)$-derivative property of intertwining operators, 
as $\bC$-formal power series with logarithmic terms, 
for $\vec{\xi}=(\theta,v,u,w)\in \Omega$ and 
$\vec{\xi}'=(\theta',v,u,w)\in \Omega'$, we have 
$$\begin{array}{l}
\frac{\pd}{\pd y}F^{A(BC)}(\vec{\xi};x,y)=F^{A(BC)}(L(-1)^{[3]}\vec{\xi};x,y),\vspace{2mm} \cr
\frac{\pd}{\pd x}F^{A(BC)}(\vec{\xi};x,y)=F^{A(BC)}(L(-1)^{[2]}\vec{\xi};x,y), \mbox{ and}\vspace{2mm}\cr
\frac{\pd}{ \pd (x-y)}F^{(AB)C}(\vec{\xi}';x-y,y)=
F^{(AB)C}(L(-1)^{[2]}\vec{\xi}';x-y,y). 
\end{array}$$
From (2.11) in Remark \ref{Derivation}, we get  
$$\begin{array}{l}
\frac{\pd}{\pd y}\{G_{h,k}^{A(BC):y}(\vec{\xi};x,y)y\}
=\frac{\pd}{\pd y}\{
F_{h,k}^{A(BC)}(\vec{\xi};x,y)y^{\gr^{234}(\vec{\xi})+1}\} \cr
\mbox{}\quad =\{\frac{\pd}{\pd y}F_{h,k}^{A(BC)}(\vec{\xi};x,y)\}y^{\gr^{234}(\vec{\xi})+1}+(\gr^{234}(\vec{\xi})+1)F_{h,k}^{A(BC)}(\vec{\xi};x,y)y^{\gr^{234}(\vec{\xi})}\cr
\mbox{}\quad =G_{h,k}^{A(BC):y}(L(-1)^{[3]}\vec{\xi};x,y)
\!-\!(k\!+\!1)G_{h,k+1}^{A(BC):y}(\vec{\xi};x,y)\}y\!+\!(\gr^{234}(\vec{\xi})\!+\!1)G_{h,k}^{A(BC):y}(\vec{\xi};x,y) \cr
\mbox{}\quad=\sum_{\vec{\xi}\in \CJ^0_N} s_{\vec{\xi}}(x,y)G_{h,k}^{A(BC):y}(\vec{\xi};x,y), 
\end{array}$$
where $s_{\vec{\xi}}(x,y)\in \bC[x,x^{-1},\iota_{x,,y}\{(x-y)^{-1}\}][y]$. 
On the other hand, we obtain:
$$
\frac{\pd}{\pd y}\{G_{h,k}^{A(BC):y}(\vec{\xi};x,y)y\}
=\frac{\pd}{\pd y}\{G_{h,k}^{A(BC):y}(\vec{\xi};x,y)\}y
+G_{h,k}^{A(BC):y}(\vec{\xi};x,y).$$
Combining it with the above, and then using Theorem 3.1, we get 
$$\begin{array}{rl}
\frac{d}{dy}\{G_{h,k}^{A(BC):y}(\vec{\xi};x_0,y)\}=&y^{-1}\frac{d}{dy}\{G^{A(BC):y}_{h,k}(\vec{\xi};x_0,y)y\}
-y^{-1}G^{A(BC):y}_{h,k}(\vec{\xi};x_0,y) \cr
=&\frac{1}{y}\sum_{\vec{\mu}\in \CJ^0_{N}}\sum_{(p,q)\in K^2}\lambda^{34:\vec{\xi},h,k}_{\vec{\mu},p,q}(x_0,y)G^{A(BC):y}_{p,q}(\vec{\mu};x_0,y)
\end{array}$$
for some $\lambda^{34:\vec{\xi},h,k}_{\vec{\mu},p,q}(x_0,y)\in \bC[x_0,x_0^{-1},
\iota_{x_0,y}\{(x_0-y)^{-1}\}][y]\subseteq \bC[[y]]$. 
Define  
$$\begin{array}{l}
G^y(x_0,y)=(G^{A(BC):y}_{h,k}(\vec{\xi};x_0,y))_{(\vec{\xi},h,k)\in \CJ_N\times K^2} 
\mbox{ and }\cr
\Lambda^{34}(x_0,y)=\left(\lambda^{34:\vec{\xi},h,k}_{\vec{\mu},p,q}(x_0,y)\right)_{(\vec{\xi},h,k),(\vec{\mu},p,q)\in \CJ_N\times K^2},\end{array}$$
then $\frac{d}{dy}G^y(x_0,y)=\frac{1}{y}\Lambda^{34}(x_0,y)G^y(x_0,y)$, 
which is a Fuchsian system over $\bD_{|x_0|}(x_0)$.
We note $\lambda^{34,\vec{\xi},h,k}_{\vec{\mu},p,q}=0$ if $\vec{\mu}\not\in \CJ^0_N$. 
Similarly, from Remark 2.3 (2.12), we have:
$$\begin{array}{l}
 (x-y_0)\frac{\pd}{\pd (x-y_0)}G^{A(BC):x-y}_{h,k}(\vec{\xi};x,y_0)=G^{A(BC):x-y}_{h,k}(L(-1)^{[2]}\vec{\xi};x,y_0)\cr
\mbox{}\qquad -(h+1)G^{A(BC):x-y}_{h+1,k}(\vec{\xi};x,y_0)x^{-1}(x-y_0)
+\gr^{234}(\vec{\xi})G^{A(BC):x-y}(\vec{\xi};x,y_0).
\end{array}\eqno{(3.2)}$$
Therefore, there are $\lambda^{23:\vec{\xi},h,k}_{\vec{\mu},p,q}(x,y_0)\in 
\bC[x,\iota_{y_0,x-y_0}\{x^{-1}\},y_0^{-1}][x-y_0]\subseteq 
\bC[[x-y_0]]$ such that 
$$\frac{\pd}{\pd x}G^{A(BC):x-y}_{h,k}(\vec{\xi};x,y_0)
=\frac{1}{x-y_0}\sum_{\vec{\mu}\in \CJ_{N}}\sum_{(p,q)\in K^2}
\lambda^{23:\vec{\xi},h,k}_{\vec{\mu},p,q}(x,y_0)G^{A(BC):x-y}_{p,q}(\vec{\mu};x,y_0). $$
Thus $G^{x-y}(x,y_0):=(G^{A(BC):x-y}_{h,k}(\vec{\xi};x,y_0))_{(\vec{\xi},p,q)\in \CJ_N\times 
K^2}$ satisfies 
a Fuchsian system 
$$ \frac{d}{dx}G^{x-y}(x,y_0)=\frac{\Lambda^{23:x-y}(x,y_0)}{x-y_0}G^{x-y}(x,y_0),$$
where $\Lambda^{23:x,y}(x,y_0)=\left(\lambda^{23:\vec{\xi},h,k}_{\vec{\mu},p,q}\right)_{(\vec{\xi},h,k),(\vec{\mu},p,q)\in \CJ_N\times K^2}$ 
and $\lambda^{23,\vec{\xi},h,k}_{\vec{\mu},p,q}=0$ if $\vec{\mu}\not\in \CJ^0_N$. 

We next investigate the set of eigenvalues of constant matrices. 
Since $\lambda^{23:\vec{\xi},h,k}_{\vec{\mu},p,q}=0$ for $\vec{\mu}\not\in \CJ_N^0$,
the set of nonzero eigenvalues of the constant matrix $\Lambda^{23;x-y}(y_0,y_0)$ for $\CJ_N$ 
is equal to that for $\CJ^0_N$. So we may assume $\CJ_N=\CJ_N^0$. 
From (3.2), we know   
$$\begin{array}{l}
 (x-y_0)\frac{\pd}{\pd (x-y_0)} G^{A(BC):x-y}_{h,k}(\vec{\xi};x,y_0)\cr
\equiv G^{A(BC):x-y}_{h,k} (L(-1)^{[2]}\vec{\xi};x,y_0)\!
+\!\gr^{234}(\vec{\xi})G_{h,k}^{A(BC):x-y}(\vec{\xi};x,y_0) \pmod{\bC[[x-y_0]](x-y_0)}. 
\end{array}$$
Furthermore, on each step of the reduction of 
$G^{A(BC):x-y}_{h,k}(L(-1)^{[2]}\vec{\xi};x,y_0)$ into a linear combination of 
$\{ G^{A(BC):x-y}_{h,k}(\vec{\mu};x,y_0)\mid \vec{\mu}\in \CJ^0_N \}$, the coefficients at the first term with ${(\alpha_j)^{\ast}}^{[1]}\vec{\xi}$ and at the 4th term with $\alpha_j^{[4]}\vec{\xi}$ are contained in 
$\bC[[x-y_0]](x-y_0)$ by Theorem 3.1. In other words, modulo $\bC[[x-y_0]](x-y_0)$,  
the process of these reductions are independent from the choice of $(\theta,w)$ 
and we don't change the subscripts $h,k$. 
Therefore, 
for $\vec{\xi}=(\theta,v,u,w), \vec{\mu}=(\tilde{\theta},\tilde{v},\tilde{u},\tilde{w})\in \CJ_N$, 
there are $\lambda^{(v,u)}_{(\tilde{v},\tilde{u})}(h,k;x,y_0)\in \bC$ such that 
$$\lambda^{23,(\theta,v,u,w),h,k}_{(\tilde{\theta},\tilde{v},\tilde{u},\tilde{w}),p,q}(x,y_0)
\equiv \delta_{h,p}\delta_{k,q}\delta_{\theta,\tilde{\theta}}\delta_{w,\tilde{w}}
\lambda^{(v,u)}_{(\tilde{v},\tilde{u})}(h,k;x,y_0) \pmod{\bC[[x-y_0]](x-y_0)}, 
\eqno{(3.3)}$$
which means that the set of eigenvalues (without counting multiplicities) of 
the constant matrix $\Lambda^{23,x-y}(y_0,y_0)$ 
is determined by only the choice of $J_{P_A}$ and $J_{P_B}$.

\subsection{Corollary 3.3 (Local normal convergence) }
Let $S^1$ be an indecomposable direct summand of $B\boxtimes C$ and $T^1$ an indecomposable direct summand of $A\boxtimes S^1$. 
We note that $A\boxtimes(B\boxtimes C)$ is a direct sum of such $T^1$'s. 
We may choose $\theta\in (T^1)^{\vee}$.  
Set $p_1=d(T^1)-d(A)-d(S^1)-1$ and $p_2=d(S^1)-d(B)-d(C)+1$. 
Then there are $\lambda^{h,k}_m(\vec{\xi})\in \bC$ such that  
$$
G^{A(BC):y}_{h,k}(\vec{\xi};x,y)= \langle \theta, \CY^1_h(v,x)\CY^2_k(u,y)w\rangle y^{\gr^{234}(\vec{\xi})}=\sum_{m=0}^{\infty} \lambda^{h,k}_m(\vec{\xi})x^{p_1}y^{p_2}(y/x)^{m},$$ 
where 
$\CY^1_h(v,z)\delta=\pi_{T^1}(\CY_h^{A(BC)}(v,z)\delta)$ and 
$\CY^2_k(u,z)w=\pi_{S^1}(\CY_k^{BC}(u,z)w)$ 
for $\delta\in S^1$ and 
$\pi_{P}$ denotes the projection to $P$.  

Therefore, in order to prove its convergence on a domain $\{(x,y)\in \bC^2 \mid 0<|y|<|x|\}$, 
it is enough to check the case where $x=1$ and $0<|y|<1$. 
In other words, it is enough to treat a one-variable vector valued function 
$Z=
\left( G^{A(BC):y}_{h,k}(\vec{\xi};1,y)\right)_{\vec{\xi}\in \CJ_{N},(h,k)\in K^2}$ satisfying a Fuchsian system $\frac{d}{dy}Z=\frac{\left(\lambda^{\vec{\xi},h,k}_{\vec{\mu},p,q}(1,y)\right)}{y}Z$ with polar locus $\{0\}$, 
where $\lambda^{\vec{\xi},h,k}_{\vec{\mu},p,q}(1,y)\in \bC[(1-y)^{-1},y]$ by Theorem 3.1. Since $(1-y)^{-1}$ and $y$ have 
radius of convergence at least $1$, so does $\lambda^{\vec{\xi},h,k}_{\vec{\mu},p,q}(1,y)-\delta_{(\vec{\xi},h,k),(\vec{\mu},p,q)}$. 
Therefore, from Theorem \ref{Radius}, the radius of convergence of $G^{A(BC):y}_{h,k}(\vec{\xi};1,y)$ (and also that of $F^{A(BC)}_{h,k}(\vec{\xi};x,y)$)
are at least $1$ for $\vec{\xi}\in \Omega$. 
In a similar way, by using $(D2)$ for 
$G^{(AB)C}(\vec{\xi}';x,y_0)$, we can show that $F^{(AB)C}(\vec{\xi}';x,y)$
is locally normally convergent on a domain $\{(x,y)\in \bC^2\mid 0<|x-y|<|y|\}$. 

This completes the proof of Corollary 3.3.

\section{Associativity} 
Set $\CD^2=\{(x,y)\in \bC^2 \mid 0<|x-y|<|y|<|x|, \mbox{ and }
x,y,x-y\not\in \bR^{\leq 0}\}$. In this section, we will prove the following main theorem. 

\begin{thm}\label{Iso} 
On $\CD^2$, we choose a principal branch $\widetilde{F}(\langle \theta,\CY^{A(BC)}(v,x)\CY^{BC}(u,y)w\rangle)$ of 
$\langle \theta,\CY^{A(BC)}(v,x)\CY^{BC}(u,y)w\rangle$ by taking the values of 
$\log(x),\log(y)$ which satisfy $-\pi <\Im(\log(x)),\Im(\log(y))<\pi$ and viewing 
$x^{d_1}y^{d_2}$ as $e^{ d_1\log(x)+d_2\log(y)}$. Similarly \\
$\widetilde{F}(\langle \theta',\CY^{(AB)C}(\CY^{AB}(v,x-y)u,y)w\rangle)$ 
is a branch of 
$\langle \theta',\CY^{(AB)C}(\CY^{AB}(v,x-y)u,y)w\rangle$ by taking the values of 
$\log(y),\log(x-y)$ which satisfy $-\pi <\Im(\log(y)),\Im(\log(x-y))<\pi$ and considering 
$y^d_1(x-y)^{d_2}=e^{d_1\log(y)+d_2\log(x-y)}$. 
Then there is an isomorphism $\phi_{[AB]C}: (A\boxtimes B)\boxtimes C
\to A\boxtimes (B\boxtimes C)$ such that \vspace{-2mm} 
$$\widetilde{F}(\langle \theta, \CY^{A(BC)}(v,x)\CY^{BC}(y,u)w\rangle) 
=\widetilde{F}(\langle \phi_{[AB]C}^{\ast}(\theta), \CY^{(AB)C}(\CY^{AB}(v,x-y)u,y)w)\rangle 
\vspace{-2mm}$$
on $\CD^2$ for any $\theta\in ((A\boxtimes B)\boxtimes C)^{\vee}$, $v\in A, u\in B, w\in C$, 
where $\phi^{\ast}_{[AB]C}$ denotes the dual of $\phi_{[AB]C}$. 
\end{thm}

\subsection{Equation (4.9)}
We may still assume $\theta\in (T^1)^{\vee}$. 
Hence there are $p_1,p_2\in \bC$ such that  
$$G^{A(BC):x-y}_{h,k}(\vec{\xi};x,y)=
\langle \theta, \CY^{A(BC)}_h(v,x)\CY^{BC}_k(u,y)w\rangle (x-y)^{\gr^{234}(\vec{\xi})}
\in \bC [[y/x]]x^{p_1}y^{p_2}$$ 
is absolutely convergent on $\{(x,y)\in \bC^2\mid 0<|y|<|x|\}$. 
We also recall the definition 
$\CJ_N$ which is the set of quadruples of basis of $
(A(BC))^{\vee}_N, \tilde{P}_A, \tilde{P}_B$, and $\tilde{P}_C$.  
We choose a branch $\widetilde{G}^{A(BC):x-y}_{h,k}(\vec{\xi};x,y)$ 
of $G^{A(BC):x-y}_{h,k}(\vec{\xi};x,y)$ 
on $\CD^2$ 
by taking values of $\log(x),\log(y)$ which satisfy 
$-\pi < \Im(\log(x)),\Im(\log(y))<\pi$ and viewing $x^{p_1+m}y^{p_2+n}$ as $e^{(p_1+m)\log(x)+(p_2+n)\log(y)}=e^{p_1\log(x)+p_2\log(y)}x^my^n$ for 
$m,n\in \bZ$. 
Then by Theorem 3.2, for $y_0\not\in \bR^{\leq 0}$, a vector valued function 
$$G^{x-y}(x,y_0)=\left(\widetilde{G}^{A(BC):x-y}_{h,k}(\vec{\xi};x,y_0)\right)
_{\vec{\xi}\in \CJ_{N}, (h,k)\in K^2}$$ 
satisfies a Fuchsian system (D2) for $\vec{\xi}\in \CJ_N$. Furthermore, 
the set of nonzero eigenvalues of the constant matrix $\Lambda^{23:x-y}(y_0,y_0)$ 
depends only on the choice of bases $J_{P_A}$ and $J_{P_B}$. 
Therefore, by Proposition \ref{Sol}, 
there is a finite subset $\Delta'\subseteq \bC$ such that all components of 
solutions of $(D2)$ are written as 
$$\sum_{t=0}^K\sum_{d\in \Delta'}\sum_{m\in \bN}r_{d,m,t}(x-y_0)^{d+m}\log^t(x-y_0) $$
with $r_{d,m,t}\in \bC$. 
We note that 
we have also taken a branch of $(x-y_0)^{d+m}\log^j(x-y_0)$ 
on $\CD^2_{(x,y_0)}$ by the same way. 
Since this holds for $\widetilde{G}^{A(BC):x-y}_{h,k}(\vec{\xi};x,y_0)$ with 
$\vec{\xi}\in \CJ_N$ and $y=y_0\not\in\bR^{\leq 0}$,  
there are $r^{h,k}_{s,t}(\vec{\xi};y)\in \bC$ and 
$K(\vec{\xi})\in \bN$ such that we can write 
$$\widetilde{G}^{A(BC):x-y}_{h,k}(\vec{\xi};x,y)=
\sum_{s\in \Delta'+\bN}\sum_{t=0}^{K(\vec{\xi})}r^{h,k}_{s,t}(\vec{\xi};y)(x-y)^{s}\log^t(x-y) \eqno{(4.1)}$$
on $\CD^2$. 
Multiplying (4.1) by $\log^h(x)\log^k(y)(x-y)^{-\gr^{234}(\vec{\xi})}$ 
for all $(h,k)\in K^2$ and taking the sum of all of them, we obtain an equality: 
$$\widetilde{F}^{A(BC)}(\vec{\xi};x,y)=
\sum_{h,k}\sum_{s\in \Delta'+\bN}\sum_{j=0}^{K(\xi)} r^{h,k}_{d,j}(\vec{\xi};y)\log^k(y)
\log^h(x)(x-y)^{s-\gr^{234}(\vec{\xi})}\log^j(x-y)  \eqno{(4.2)}$$
for a branch $\widetilde{F}^{A(BC)}(\vec{\xi};x,y)$ of $F^{A(BC)}(\vec{\xi};x,y)$ on $\CD^2$. 
About the derivations, we use $\frac{\pd x}{\pd(x-y)}=1$ and $\frac{\pd y}{\pd(x-y)}=0$. 
Since $\frac{\pd \log(x)}{\pd x}=\frac{1}{x}=\frac{1}{y+(x-y)}=\frac{1}{y}\sum_{j=0}^{\infty}(-1)^j(\frac{x-y}{y})^j$ and 
$\frac{\pd \log(x)}{\pd y}=0$, we define  
$$ p(x,y):=\sum_{j=0}^{\infty}\frac{(-1)^j}{j+1}\left(\frac{x-y}{y}\right)^{j+1}+\log(y), \eqno{(4.3)}$$
which satisfies the same property with $\log(x)$ in a neighborhood of $y\not\in \bR^{\leq 0}$. We replace 
$\log(x)$ in (4.2) by $p(x,y)$. Since the powers of $(x-y)$ in $p(x,y)$ are non-negative, 
we can write 
$$ \widetilde{F}^{A(BC)}(\vec{\xi};x,y)
=\sum_{t=0}^{K(\vec{\xi})} \sum_{s\in \gr^{234}(\vec{\xi})-\Delta'-\bN-1}g_{s,t}(\vec{\xi};y)(x-y)^{-s-1}\log^t(x-y)  
\eqno{(4.4)}$$
with $g_{s,t}(\vec{\xi};y)\in \bC$ for $\vec{\xi}\in \CJ_N$.  
Since $\dim P_C<\infty$, 
using $\Delta=\{ d-{\rm max}\{\gr(w):w\in P_C\}\mid d\in \Delta'\}$, 
we can rewrite (4.4) into 
$$ \widetilde{F}^{A(BC)}(\vec{\xi};x,y)
=\sum_{t=0}^{K(\vec{\xi})} \sum_{s\in \gr(v)+\gr(u)-\Delta-\bN-1}g_{s,t}(\vec{\xi};y)(x-y)^{-s-1}\log^t(x-y)  
\eqno{(4.5)}$$
for $\vec{\xi}=(\theta,v,u,w)\in \Omega$ with $w\in P_C$.

\begin{lmm}
$(4.5)$ holds for all $\vec{\xi}\in \Omega$. 
\end{lmm}

\pr 
Suppose false and let $\vec{\xi}=(\theta,v,u,w')$ be a counterexample. 
Subject to counterexample, we choose $w'$ with the smallest grade. 
Since $(4.5)$ holds for elements with $w$ in $P_C$, we may assume $w'=\alpha_{-1}w$ with $\alpha\in V_{\geq 0}$ and 
$g_{s,t}(\theta,v,u,\alpha_{-1}w;y)\not=0$ for some $s\not\in \wt(v)+\wt(u)-1-\Delta-\bN$. In this case, by Borcherds identities, we obtain  
$$\begin{array}{l}
\widetilde{F}^{A(BC)}(\alpha_{-1}^{[4]}\vec{\xi};x,y)=
\widetilde{F}^{A(BC)}({(\alpha_{-1})^{\ast}}^{[1]}\vec{\xi};x,y) 
\!-\!\sum_{j=0}^{\infty}(-1)^j\widetilde{F}^{A(BC)}(\alpha_j^{[2]}\vec{\xi};x,y)x^{-j-1}\cr
\mbox{}\,\, -\sum_{j=0}^{\infty}(-1)^j
\widetilde{F}^{A(BC)}(\alpha_j^{[3]}\vec{\xi};x,y)y^{-j-1}\cr 
=\widetilde{F}^{A(BC)}({(\alpha_{-1})^{\ast}}^{[1]}\vec{\xi};x,y)-\sum_{j=0}^{\infty}\binom{-1}{j}\widetilde{F}^{A(BC)}(\alpha_j^{[2]}\vec{\xi};x,y)x^{-j-1}\cr
\mbox{}\,\, -\sum_{j=0}^{\infty}\binom{-1}{j}y^{-j-1}\{\sum_{i\in \bN}\binom{j}{i}(-1)^i\widetilde{F}^{A(BC)}
({(\alpha_{j-i}^{\ast}}^{[1]}\vec{\xi};x,y)y^i \cr
\mbox{}\,\, -\sum_{i\in \bN}\binom{j}{i}
\widetilde{F}^{A(BC)}(\alpha_i^{[2]}\vec{\xi};x,y)(x-y)^{j-i} 
\!-\!\sum_{i\in \bN}\binom{j}{i}(-1)^{i+j+1}\widetilde{F}^{A(BC)}(\alpha_i^{[4]}\vec{\xi};x,y)y^{j-i}\} \cr
=\widetilde{F}^{A(BC)}({(\alpha_{-1})^{\ast}}^{[1]}\vec{\xi};x,y)
-\sum_{j,i\in \bN}\binom{j}{i}(-1)^{i+j}
y^{i-j-1}\widetilde{F}^{A(BC)}({(\alpha_{j-i})^{\ast}}^{[1]}\vec{\xi};x,y)\cr
\mbox{}\,\,
+\sum_{j=0}^{\infty}\sum_{i\in \bN}\binom{j}{i}(-1)^{i+1}\widetilde{F}^{A(BC)}
(\alpha_i^{[4]}\vec{\xi};x,y)y^{-1}, 
 \end{array}$$
by (2.9). However, since $\gr(\alpha_jw)<\gr(\alpha_{-1}w)=\gr(w')$ for $j\geq 0$ and 
the 2nd and 3rd components of $(\alpha_j)^{[i]}\vec{\xi}$ in RHS of 
the last equation are $v$ and $u$, none of the terms in RHS of the last equation have nonzero coefficients of $(x-y)^{-s-1}\log^t(x-y)$ 
by the minimality of $\gr(w')$, which is a contradiction. \qed

Since $g_{s,t}(\vec{\xi};y)$ is multi-linear on 
$\vec{\xi}\in \Omega_N=(A(BC))^{\vee}_{(\leq N)}\times A\times B\times C$, 
we split $\vec{\xi}=(\theta,v,u,w)$ into $\{v,u\}$ and $\{\theta, w\}$. In other words,
using a formal operator $\widetilde{\CY}^{3}_{s,t}\in \Hom(A\otimes B, \Hom(C, A(BC)\otimes 
O(\bC\backslash\bR^{\leq 0})))$, we can write 
$$\widetilde{F}^{A(BC)}(\vec{\xi};x,y)
=\sum_{t=0}^{K(\vec{\xi})}\sum_{s\in \gr(v)+\gr(u)-\Delta-\bN-1}\langle \theta, 
\widetilde{\CY}^{3}_{s,t}(v,u;y)w\rangle (x-y)^{-s-1}\log^t(x-y)  \eqno{(4.6)}$$ 
for each $(x,y)\in \CD^2$. Namely, 
$\langle \theta,\widetilde{\CY}^3_{s,t}(v,u;y)w\rangle=\Coeff_{s,t}(\tilde{F}^{A(BC)}(\vec{\xi};x,y);y)$, where \\$\Coeff_{s,t}(f(x,y);y)$ denotes a coefficient $\lambda_{s,t}(y)$ of $(x-y)^{-s-1}\log^t(x-y)$ in \\$f(x,y)=
\sum_{t\in \bN}\sum_{s\in \bC} \lambda_{s,t}(y)(x-y)^{-s-1}\log^t(x-y)$.  

We introduce a formal vector space $A\otimes_{s,t}B$, which is isomorphic to $A\otimes B$ as a vector space and its isomorphism is given by $v\otimes v\mapsto v\otimes_{s,t}u$ for $s,t$. 
Set 
$A\otimes_{\infty}B=\oplus_{(s,t)\in \bC\times \bN}(A\otimes_{s,t}B)$ and define  
a formal operator $\widetilde{\CY}^4$ of type $\binom{A\otimes_{\infty}B}{A\,\,\,\, B}$ 
by $\widetilde{\CY}^4(v,x-y)u=\sum_{s\in \bC}\sum_{t\in \bN}v\otimes_{s,t}u (x-y)^{-s-1}\log^t(x-y)$. 
We also define another formal operator 
$\widetilde{\CY}^3: (A\otimes_{\infty}B)\times C\to A(BC)\otimes \CO(\CD^2_{(x,y)})$ 
by $\langle \theta, \widetilde{\CY}^3(v\otimes_{s,t}u,y)w\rangle
=\langle \theta, \widetilde{\CY}^3_{s,t}(v,u;y)w\rangle (x-y)^{-s-1}\log^t(x-y)$.  
Then we can rewrite Eq. (4.6) into  
$$
\widetilde{F}^{A(BC)}(\vec{\xi};x,y)=\langle \theta, 
\widetilde{\CY}^3(\widetilde{\CY}^4(v,x-y)u,y)w\rangle. \eqno{(4.7)}$$

As it is well-known, $V[z,z^{-1}]$ has a Lie algebra structure by \vspace{-2mm}
$$[\alpha z^m, 
\beta z^n]=\sum_{j\in \bN} \binom{m}{j}(\alpha_j\beta) z^{n+m-j}.$$

\begin{lmm} 
Define the action of $V[z,z^{-1}]$ on $A\otimes_{\infty}B$ by  
$(\alpha z^m)(v\otimes_{s,t}u)=v\otimes_{s,t}\alpha_mu+\sum_{j=0}^{\infty}\binom{m}{j}
(\alpha_jv\otimes_{s+m-j,t} u) \mbox{ for }\alpha z^m \in V[z,z^{-1}]$.
Then $A\otimes_{\infty}B$ is a $V[z,z^{-1}]$-module. 
\end{lmm}

\pr  
To simplify the notation, we will omit the index $t$ from $v\otimes_{s,t}u$. 
We also denote the action of 
$\alpha z^m$ by $\alpha_{[m]}$.  
Then for $\alpha,\beta\in V$, we obtain 
$$\begin{array}{l}
 \alpha_{[m]}\beta_{[n]}(v\otimes_su)=\alpha_{[m]}(\sum_{j=0}^{\infty}\binom{n}{j}(\beta_jv)\otimes_{n+s-j}u)+\alpha_{[m]}(v\otimes_s(\beta_nu)) \cr
\mbox{}\quad =\sum_{i=0}^{\infty}\sum_{j=0}^{\infty}\binom{m}{i}\binom{n}{j}(\alpha_i\beta_jv)\otimes_{m+n+s-i-j}u
 +\sum_{j=0}^{\infty}\binom{n}{j}(\beta_jv)\otimes_{n+s-j}(\alpha_mu)\cr
\mbox{}\qquad+\sum_{i=0}^{\infty}\binom{m}{i}(\alpha_iv)\otimes_{m+s-i}(\beta_nu)+v\otimes_s(\alpha_m\beta_nu).
\end{array}$$
We hence get:
$$\begin{array}{l}
 [\alpha_{[m]},\beta_{[n]}](v\otimes_su)\cr
\mbox{}\quad=\sum_{i=0}^{\infty}\sum_{j=0}^{\infty}\binom{m}{i}\binom{n}{j}(\alpha_i\beta_jv)\otimes_{m+n+s-i-j}u
 +\sum_{j=0}^{\infty}\binom{n}{j}(\beta_jv)\otimes_{n+s-j}(\alpha_mu)\cr
\mbox{}\qquad+\sum_{i=0}^{\infty}\binom{m}{i}(\alpha_iv)\otimes_{m+s-i}(\beta_nu)+v\otimes_s(\alpha_m\beta_nu)\cr
\mbox{}\qquad-\sum_{j=0}^{\infty}\sum_{i=0}^{\infty}\binom{n}{j}\binom{m}{i}(\beta_j\alpha_iv)\otimes_{m+n+s-i-j}u
 -\sum_{i=0}^{\infty}\binom{m}{i}(\alpha_iv)\otimes_{m+s-i}(\beta_nu)\cr
\mbox{}\qquad-\sum_{j=0}^{\infty}\binom{n}{j}(\beta_jv)\otimes_{n+s-j}(\alpha_mu)-v\otimes_s(\beta_n\alpha_mu) \cr
\mbox{}\quad=\sum_{i=0}^{\infty}\sum_{j=0}^{\infty}\binom{m}{i}\binom{n}{j}([\alpha_i,\beta_j]v)\otimes_{m+n+s-i-j}u
+v\otimes_s([\alpha_m,\beta_n]u) \cr
\mbox{}\quad=\sum_{i=0}^{\infty}\sum_{j=0}^{\infty}\binom{m}{i}\binom{n}{j}([\alpha_i,\beta_j]v)\otimes_{m+n+s-i-j}u
+v\otimes_s([\alpha_m,\beta_n]u) \cr
\mbox{}\quad=\sum_{i=0}^{\infty}\sum_{j=0}^{\infty}\binom{m}{i}\binom{n}{j}(\sum_{p=0}^{\infty}\binom{i}{p}(\alpha_p\beta)_{i+j-p}v)\otimes_{m+n+s-i-j}u
+v\otimes_s([\alpha_m,\beta_n]u).
\end{array}\eqno{(4.8)} $$
On the other hand, we have:
$$\begin{array}{l}
\sum_{p=0}^{\infty}\binom{m}{p}(\alpha_p\beta)_{[m+n-p]}(v\otimes_su) \cr
\mbox{}\quad =v\otimes_s(\sum_{p=0}^{\infty}\binom{m}{p}(\alpha_p\beta)_{m+n-p}u)
+\sum_{p=0}^{\infty}\binom{m}{p}\sum_{q=0}^{\infty}\binom{m+n-p}{q}((\alpha_p\beta)_qv)
\otimes_{s+m+n-p-q}u \cr
\mbox{}\quad =\sum_{i=0}^{\infty}\sum_{j=0}^{\infty}\binom{m}{i}\binom{n}{j}(\sum_{p=0}^{\infty}\binom{i}{p}(\alpha_p\beta)_{i+j-p}v)\otimes_{m+n+s-i-j}u
+[v\otimes_s([\alpha_m,\beta_n]u),  
\end{array}$$
which coincides with the last terms in (4.8). 
Therefore $A\otimes_{\infty}B$ is 
a $V[z,z^{-1}]$-module. 
\prend

$A\otimes_{\infty}B$ is not a $V$-module. We choose a better one.

\begin{defn}
Let $\Gamma$ be the kernel of $\widetilde{\CY}^3$, that is, 
$$\Gamma=\{ \gamma\in A\otimes_{\infty}B \mid  
\langle \theta, \widetilde{\CY}^3(\gamma, y)w\rangle=0 
\mbox{ for all }\theta\in (A(BC))^{\vee}, w\in C \mbox{ and }y\not\in \bR^{\leq 0}\}.$$ 
In other words, $\sum_{i=1}^n v^i\otimes_{s_i,t_i}u^i\in \Gamma$ if and only if  
$\sum_{i=1}^n \Coeff_{s_i,t_i}(\widetilde{F}^{A(BC)}(\theta, v^i,u^i,w;x,y);y)=0$ 
for any $\theta\in (A(BC))^{\vee}$ and $w\in C$. 
Clearly, $\Gamma$ is a subspace of $A\otimes_{\infty}B$. 
Define 
$A\tilde{\boxtimes} B=(A\otimes_{\infty}B)/\Gamma$. 
We denote $v\otimes_{s,t}u+\Gamma$ 
in $A\tilde{\boxtimes}B$ by $[v\otimes_{s,t}u]$ or $[v_{s,t}u]$. 
\end{defn}

\begin{lmm} $\widetilde{\CY}^3$ satisfies associativity ${\rm (I\,\,2)}$ and 
$\Gamma$ is $V[z,z^{-1}]$-invariant.  
\end{lmm}

\pr 
By the definition of the action of $V[z,z^{-1}]$, we get:
$$\begin{array}{l}
\langle \theta, \widetilde{\CY}^3(\alpha_{[n]}(v\otimes_{s,t}u),y)w\rangle \cr
\mbox{}\quad=\langle \theta, \widetilde{\CY}^3(v\otimes_{s,t}\alpha_n u),y)w\rangle
+\sum_{i=0}^{\infty}\binom{n}{i}\langle \theta, 
\widetilde{\CY}^3(\alpha_iv\otimes_{s+n-i,t}u,y)w\rangle \cr
\mbox{}\quad=\Coeff_{s,t}(\widetilde{F}^{A(BC)}(\theta,v,\alpha_nu,w;x,y);y)\!+\!\sum_{i=0}^{\infty}
\binom{n}{i}\Coeff_{s+n-i,t}(\widetilde{F}^{A(BC)}(\theta,\alpha_jv,u,w;x,y);y).
\end{array}$$
Since $\widetilde{F}^{A(BC)}$ satisfies Borcherds identity (3B), we also have:
$$\begin{array}{l}
\Coeff_{s,t}(\widetilde{F}^{A(BC)}(\theta,v,\alpha_nu,w;x,y);y)=\Coeff_{s,t}\left(\left( \sum_{j=0}^{\infty}(-1)^j\widetilde{F}^{A(BC)}(((\alpha_{n-j})^{\ast})^{[1]}\vec{\xi};x,y)y^j \right.\right.\cr
\left.\left.\mbox{}\quad\!-\!\sum_{j=0}^{\infty}\binom{n}{j}\widetilde{F}^{A(BC)}(\alpha_j^{[2]}\vec{\xi};x,y)(x\!-\!y)^{n-j}\!-\!\sum_{j=0}^{\infty}\binom{n}{j}\widetilde{F}^{A(BC)}(
\alpha_j^{[4]}\vec{\xi};x,y)(-y)^{n-j}\right);y\right)\cr
\mbox{} =\sum_{j=0}^{\infty}(-1)^j
\Coeff_{s,t}(\widetilde{F}^{A(BC)}(((\alpha_{n-j})^{\ast})^{[1]}\vec{\xi};x,y);y)y^j\cr
\mbox{}\quad
-\sum_{j=0}^{\infty}\binom{n}{j}\Coeff_{s+n-j,t}(\widetilde{F}^{A(BC)}(
\alpha_j^{[2]}\vec{\xi};x,y);y)\cr
\mbox{}\quad
-\sum_{j=0}^{\infty}\binom{n}{j}(-1)^{j+n}\Coeff_{s,t}(\widetilde{F}^{A(BC)}(\alpha_j^{[4]}\vec{\xi};x,y);y)y^{n-j}. 
\end{array}$$
Therefore, we obtain ${\rm (I\,\,2)}$ for $\widetilde{\CY}^3$, that is, 
$$
\begin{array}{l}
\langle \theta, \widetilde{\CY}^3(\alpha_{[n]}(v\otimes_{s,t}u), y)w\rangle
=\sum_{j=0}^{\infty}\binom{n}{j} (-1)^j\langle (\alpha_{n-j})^{\ast}\theta, 
\widetilde{\CY}^3(v\otimes_{s,t}u, y)w\rangle y^j \cr
\mbox{}\qquad -\sum_{j=0}^{\infty} \binom{n}{j}(-1)^{j+n}
\langle \theta,\widetilde{\CY}^3(v\otimes_{s,t}u,y)\alpha_j w \rangle y^{n-j}.
\end{array}$$
We next show that $\Gamma$ is $V[z,z^{-1}]$-invariant. 
If $\gamma=\sum_{i=1}^r v^i\otimes_{s_i,t_i}u^i\in \Gamma$, then we have  
$$\begin{array}{l}
\langle \theta, \widetilde{\CY}^3(\alpha_n\gamma,y)w\rangle \cr
=\sum_{j=0}^{\infty}\binom{n}{j}(-1)^j 
\langle (\alpha_{n-j})^{\ast}\theta, \widetilde{\CY}^3(\gamma,y)w\rangle y^j
-\sum_{j=0}^{\infty}\binom{n}{j}(-1)^{j+n}\langle \theta, 
\widetilde{\CY}^3(\gamma,y)\alpha_j w\rangle y^{n-j}=0
\end{array}$$ 
for any $\theta$ and $w$, which means 
$\alpha_n\gamma\in \Gamma$ for any $\alpha z^n=\alpha_n$. Hence $\Gamma$ is $V[z,z^{-1}]$-invariant.  
\prend

Clearly, we can define formal operators  
$\CY^{3}$ of type $\binom{A(BC)}{A\tilde{\boxtimes}B\,\, C}$ and 
$\CY^{4}$ of type $\binom{A\tilde{\boxtimes}B}{A\,\, B}$ by  
$\CY^{4}(v,z)u=[\widetilde{\CY}^4(v,z)u]=\sum_{s,t}[v\otimes_{s,t}u]z^{-s-1}\log^t(z)$ and 
$\langle \theta, \CY^{3}([v_{s,t}u],z)w\rangle=\langle \theta, 
\widetilde{\CY}^3(v\otimes_{s,t}u,z)w\rangle$ for $\vec{\xi}\in \Omega$. 
Then we have 
$$ \widetilde{F}^{A(BC)}(\vec{\xi};x,y)
=\langle \theta, \CY^{3}(\CY^{4}(v,x-y)u,y)w\rangle  \eqno{(4.9)}$$
for $\vec{\xi}\in \Omega$. We use $F^{34}(\vec{\xi};x,y)$ to denote 
$\langle \theta, \CY^{3}(\CY^{4}(v,x-y)u,y)w\rangle$.

\begin{thm}
$\langle \theta, \CY^3(\CY^4(v,x-y)u,y)w\rangle$ 
satisfies the Borcherds identities $(1B)\sim(4B)$. 
\end{thm}

\pr 
As we mentioned in Comment 2, the coefficients of 
$F^{A(BC)}(\vec{\mu};x,y)$ in the expansion of 
$F^{A(BC)}(\alpha^{[i]}_a\vec{\xi};x,y)$ by Borcherds identities 
are all $\bC[x^{\pm 1},y^{\pm 1}, 
(x-y)^{\pm 1}]$, that is, integer powers of $x, y$ and $x\!-\!y$. 
Hence branches $\widetilde{F}^{A(BC)}(\vec{\xi};x,y)$ of them on $\CD^2$ 
also satisfy the Borcherds identities $(1A)\sim (4A)$. 
By Eq. (4.9), $F^{34}(\vec{\xi};x,y)$ also satisfies 
the Borcherds identities $(1A)\sim (4A)$ on $\CD^2$.  
Since the corresponding coefficients in $(1A)\sim (4A)$ and $(1B)\sim (4B)$ are same on $\CD^2$, we have the desired identities. \qed

\begin{prop}
$\CY^{3}$ and $\CY^{4}$ satisfy $L(-1)$-derivative properties. 
\end{prop}
 
\pr
We note that we are considering the derivations satisfying 
$\frac{\pd}{\pd (x-y)}y=0=\frac{\pd}{\pd y}(x-y)$ and 
$\frac{\pd}{\pd y}y=1=\frac{\pd}{\pd(x-y)}(x-y)$. 
By the choice of $p(x,y)$ in (4.4), we get  
$$
\begin{array}{l}
\langle \theta, 
\CY^{3}(\CY^{4}(L(-1)v,x-y)u,y)w\rangle
=\widetilde{F}^{A(BC)}(L(-1)^{[2]}\vec{\xi};x,y) 
=\frac{d}{d x}\widetilde{F}^{A(BC)}(\vec{\xi};x,y) \cr
\mbox{}\qquad=\frac{d}{d (x-y)}\widetilde{F}^{A(BC)}(\vec{\xi};x,y)
=\frac{\pd }{\pd (x-y)}\langle \theta,\CY^{3}(\CY^{4}(v,x-y)u,y)w\rangle.   
\end{array}$$
We fix $v\in A$. 
Since the above holds 
for every $\theta\in (A(BC))^{\vee}, u\in B, w\in C$, the above implies 
$L(-1)$-derivative property of $\CY^4$: 
$\CY^4(L(-1)v,x-y)=\frac{d}{d(x-y)} \CY^4(v,x-y)$. Using \\
$F^{34}(L(-1)^{[3]}\vec{\xi};x,y)=\frac{d}{d y} \widetilde{F}^{A(BC)}(\vec{\xi};x,y)
=\frac{\pd }{\pd y}F^{34}(\vec{\xi};x,y)$, we obtain 
$$\begin{array}{l}
\langle \theta, 
\CY^{3}(L(-1)\CY^{4}(v,x-y)u,y)w\rangle\cr
\mbox{}\qquad =\langle \theta, \CY^{3}(\CY^{4}(L(-1)v,x-y)u,y)w\rangle
+\langle \theta,\CY^{3}(\CY^{4}(v,x-y)L(-1)u,y)w\rangle \cr
\mbox{}\qquad =\frac{\pd}{\pd x}\widetilde{F}^{A(BC)}(\vec{\xi};x,y)
+\frac{\pd}{\pd y}\widetilde{F}^{A(BC)}(\vec{\xi};x,y)=(\frac{\pd}{\pd x}+\frac{\pd}{\pd y})\widetilde{F}^{A(BC)}(\vec{\xi};x,y). 
\end{array}$$
Since $(\frac{\pd}{\pd y}+\frac{\pd}{\pd x})(x-y)=0$ and 
$(\frac{\pd}{\pd y}+\frac{\pd}{\pd x})y=1$, we also get 
$$\begin{array}{l}
\langle \theta, 
\CY^{3}(L(-1)\CY^{4}(v,x-y)u,y)w\rangle
=\frac{\pd}{\pd y}\CY^{3}(\CY^{4}(v,x-y)u,y)w\rangle
\end{array}$$
under the assumption $\frac{\pd (x-y)}{\pd y}=0$, which proves 
$L(-1)$-derivative property of $\CY^{3}$. 
\prend

As a corollary of the above and the existence of a finite set $\Delta$ in (4.5), we have: 

\begin{cry}
$L(-1)[v_{s,t}u]=-s[v_{s-1,t}u]+(t+1)[v_{s-1,t+1}u]$ and 
$\wt([v_{s,t}u])=\wt(v)+\wt(u)-s-1$. The weights of elements in $A\tilde{\boxtimes}B$ are contained in $\Delta+\bN$ and so $A\tilde{\boxtimes}B$ is $\bN$-gradable. 
In particular, $\alpha_n[v_{s,t}u]=0$ for $n>\!>0$ for $\alpha\in V$. 
\end{cry}

Although Borcherds identities are given from ${\rm (I\,\,1)}$ and ${\rm (I\,\, 2)}$, 
we will show the reverse.  

\begin{prop}
$\CY^{3}$ and $\CY^{4}$ satisfy ${\rm (I\,\, 1)}$ and ${\rm (I\,\, 2)}$. 
\end{prop}

\pr 
We have already proved ${\rm (I\,\, 2)}$ for $\CY^3$. 
By the action of $V$ on $A\otimes_{\infty}B$, we have: 
$$\begin{array}{l}\langle \theta, \CY^{3}(\alpha_{[n]}\CY^{4}(v,x-y)u,y)w\rangle
=\langle \theta, \widetilde{\CY}^3(\alpha_n\widetilde{\CY}^4(v,x-y)u,y)w\rangle \cr
\mbox{}\quad=\langle \theta, \widetilde{\CY}^3(\widetilde{\CY}^4(v,x-y)\alpha_nu,y)w\rangle 
+\sum_{j=0}^{\infty}\binom{n}{j}\langle \theta, \widetilde{\CY}^3(\widetilde{\CY}^4(\alpha_jv,x-y)u,y)w\rangle (x-y)^{n-j}  
\end{array}$$
for all $\theta$ and $w$, which implies ${\rm (I\,\,1)}$ for $\CY^4$. 
From Borcherds identity (2B), we obtain:
$$\begin{array}{l}
\langle \theta, \CY^{3}(\CY^{4}(\alpha_nv,x-y)u,y)w\rangle
=\sum_{j=0}^{\infty}\binom{n}{j}(-1)^j\langle \theta, \alpha_{n-j}\CY^{3}(\CY^{4}(v,x-y)u,y)w\rangle x^j\cr
\mbox{}\quad -\sum_{j=0}^{\infty}\binom{n}{j}(-1)^{n+j}\langle\theta,\CY^{3}(\CY^{4}(v,x-y)u,y)\alpha_jw\rangle x^{n-j} \cr
\mbox{}\quad -\sum_{j=0}^{\infty} \binom{n}{j}(-1)^{n+j}\langle \theta, \CY^{3}(\CY^{4}(v,x-y)\alpha_ju,y)w\rangle(x-y)^{n-j}. 
\end{array}$$
On the other hand, from the associativity ${\rm (I\,\,2)}$ for $\CY^{3}$, we have 
$$\begin{array}{l}
\sum_{j=0}^{\infty}\binom{n}{j}(-1)^j\langle \theta, \CY^{3}(\alpha_{n-j}\CY^{4}(v,x-y)u,y)w\rangle(x-y)^j \cr 
\mbox{}\quad =\sum_{i,j\in \bN}\binom{n}{j}\binom{n-j}{i}(-1)^{i+j}\langle \theta, \alpha_{n-i-j}\CY^{3}(\CY^{4}(v,x-y)u,y)w\rangle(x-y)^jy^i \cr
\mbox{}\qquad -\sum_{i,j\in \bN}\binom{n}{j}(-1)^{i+n}\binom{n-j}{i}\langle \theta, 
\CY^{3}(\CY^{4}(v,x-y)u,y)\alpha_jw\rangle (x-y)^jy^{n-j-i}.
\end{array}$$
Since $\sum_{j=0}^{\ell}\binom{n}{j}\binom{n-j}{\ell-j}(x-y)^jy^{\ell-j}=\binom{n}{\ell}x^{\ell}$ for $\ell\in \bN$ and 
$\sum_{j}\binom{n}{j}\binom{n-j}{i}(-1)^{i+n}(x-y)^jy^{n-j-i}=\binom{n}{i}(-1)^{n-i}x^{n^i}$ by (2.9), 
we obtain ${\rm (I\,\, 2)}$ for $\CY^{4}$. 
Commutativity ${\rm (I\,\, 1)}$ for $\CY^{3}$ is a direct consequence of Borcherds identity 
$(B1)$ and the definition of action of $V$ on $A\tilde{\boxtimes}B$. \qed

\subsection{Recovery of $V$-modules and intertwining operators}
\begin{thm}\label{Vmod}
$A\tilde{\boxtimes}B$ is an $\bN$-gradable $V$-module. 
\end{thm}

\pr  
We have proved that $A\tilde{\boxtimes}B$ is an $\bN$-gradable $V[z,z^{-1}]$-module in Corollary 4.8. In particular, the action of $V$ on $A\tilde{\boxtimes} B$ satisfies commutativity. 
Therefore, the remaining thing is to prove associativity of the action of $V$. 
In particular, it is enough to show 
$$\begin{array}{l}
\langle \theta, \CY^3(\alpha_{-1}\beta)_{[-1]}\gamma,y)w\rangle \cr
\mbox{}\quad =\langle \theta, \CY^3(\alpha_{[-1]}\beta_{[-1]}\gamma,y)w\rangle
 +\sum_{j=0}^{\infty}\langle \theta, (\alpha_{[-2-j]}\beta_{[j]}+\beta_{[-2-j]}\alpha_{[j]})\gamma,y)w\rangle
\end{array} \eqno{(4.10)}$$
for $\gamma\in A\tilde{\boxtimes}B, w\in C, \theta\in (A(BC))^{\vee}$. 
We will develop the action of $V$ into a triple $(A(BC))^{\vee}\times A\tilde{\boxtimes}B \times C$ by using ${\rm (I\,\,1)}$ and ${\rm (I\,\,2)}$. 
The rule is simple, that is, we shift $\alpha_{<0}$ and $\beta_{<0}$ to the left and $\alpha_{\geq 0}$ and 
$\beta_{\geq 0}$ to the right. Finally, the actions of $V$ on $A\tilde{\boxtimes}B$ and 
$C$ are given by $\alpha_{\geq 0}$ or $\beta_{\geq 0}$. 
We use notation $P(X(\alpha,\beta))$ to denote $X(\alpha,\beta)+X(\beta,\alpha)$ 
for any form $X$ of $\alpha_i$ and $\beta_j$. 
For example, $P(\langle \theta, \CY^3(\alpha_{[-2-j]}\beta_{[j]})\gamma,y)w\rangle)$ denotes 
$\langle \theta, \CY^3(\alpha_{[-2-j]}\beta_{[j]}
+\beta_{[-2-j]}\alpha_{[j]})\gamma,y)w\rangle$.
Then by ${\rm (I\,\,2)}$ and associativity of the action of $V$ on $A(BC)$ and $C$, we have:
$$\begin{array}{l}
{\rm LHS}\mbox{ of (4.10)}=\sum_{j\in \bN}\langle\theta, (\alpha_{-1}\beta)_{-1-j}\CY^3(\gamma,y)w\rangle y^j
+\sum_{j\in \bN} \langle \theta, \CY^3(\gamma,y)(\alpha_{-1}\beta)_jw\rangle y^{-1-j}\cr
\mbox{}=\sum_{i,j\in \bN}\langle \theta, 
\{ \alpha_{-1-j}\beta_{-1-j+i}+\beta_{-2-j-i}\alpha_i\}\CY^3(\gamma,y)w\rangle y^j \cr
\mbox{}\,\, +\sum_{i,j\in \bN} 
\langle \theta, \CY^3(\gamma,y)\{\alpha_{-1-i}\beta_{j+i}+\beta_{-1+j-i}\alpha_i\}w\rangle y^{-1-j}\cr
\mbox{}=\sum_{\ell,i\in \bN}\langle \theta, \alpha_{-1-i}\beta_{-1-\ell}\CY^3(\gamma,y)w\rangle y^{\ell+i}
+\sum_{\ell,j\in \bN}\langle \theta, \alpha_{-2-j-\ell} \beta_{\ell}\CY^3(\gamma,y)w\rangle y^j\cr
\mbox{}\,\,+\sum_{i,j\in \bN}\langle \theta, \beta_{-2-i-j}\alpha_i\CY^3(\gamma,y)w\rangle y^j  
+\sum_{i,j\in \bN}\langle \theta, \CY^3(\gamma,y)\alpha_{-1-i}\beta_{i+j}w\rangle y^{-1-j}\cr
\mbox{}\,\, +\sum_{i,j\in \bN}\langle \theta, \CY^3(\gamma,y)
\beta_{-1-\ell}\alpha_{\ell+j}w\rangle y^{-1-j}+\sum_{\ell,i\in \bN}\langle \theta, \CY^3(\gamma,y)\beta_{\ell}\alpha_iw\rangle y^{-2-\ell-i} \cr
\mbox{}=\sum_{\ell,i\in \bN}\langle \theta, \alpha_{-1-i}\beta_{-1-\ell}\CY^3(\gamma,y)w\rangle y^{\ell+i}
+P(\sum_{i,j\in \bN}\langle \theta, \beta_{-2-i-j}[\alpha_i,\CY^3(\gamma,y)]w\rangle y^j) \cr
\mbox{}\,\, +\!P(\sum_{i,j\in \bN}\langle \theta, \beta_{-2-i-j}\CY^3(\gamma,y)\alpha_iw\rangle y^j)\!
-\!P(\sum_{i,j\in \bN}\langle \theta, [\alpha_{-1-i},\CY^3(\gamma,y)]\beta_{i+j}w\rangle y^{-1-j}) \cr
\mbox{}\,\,+\!P(\sum_{i,j\in \bN}\langle \theta, \alpha_{-1-i}\CY^3(\gamma,y)\beta_{i+j}w\rangle y^{-1-j})
+\sum_{\ell,i\in \bN}\langle \theta, \CY^3(\gamma,y)\beta_{\ell}\alpha_iw\rangle y^{-2-\ell-i}. 
\end{array}$$
Since there are many subscripts, in order to count the coefficients of terms in summations, we use monomials. For example, we denote 
$\langle \theta, \beta_{a}\CY^3([v_su],y)\alpha_b w\rangle y^c$
by $X^{a}Y^bZ^c$. 
Then since $i,j\in \bN$, keeping the non-negative powers of $Y$, 
the following equation
$$\begin{array}{l}
\sum_{i,j\in \bN}X^{-2-i-j}Y^iZ^j+\sum_{i,j\in \bN}X^{-1-i}Y^{i+j}Z^{-1-j} \cr
\mbox{}\quad=\iota_{X,Y}\{(X-Y)^{-1}\}\iota_{X,Z}\{(X-Z)^{-1}\}+\iota_{X,Y}\{(X-Y)^{-1}\}\iota_{Z,Y}\{(Z-Y)^{-1}\}\cr
\mbox{}\quad =\iota_{X,Z}\{(X-Z)^{-1}\}\iota_{X,Y}\{(X-Y)^{-1}\}\iota_{Z,Y}\{(Z-Y)^{-1}\}(X-Y)\cr
\mbox{}\quad =\iota_{X,Z}\{(X-Z)^{-1}\}\iota_{Z,Y}\{(Z-Y)^{-1}\}\cr
\mbox{}\quad =\sum_{i,j\in \bN}X^{-1-i}Y^jZ^{i-j-1} 
\end{array}$$
implies 
$$\begin{array}{l} 
\sum_{i,j\in \bN}\langle \theta, \beta_{-2-i-j}\CY^3(\gamma,y)\alpha_iw\rangle y^j
+\sum_{i,j\in \bN}\langle \theta, \beta_{-1-i}\CY^3(\gamma,y)\alpha_{i+j}w\rangle y^{-1-j}\cr
\mbox{}\qquad =\sum_{i,j\in \bN}\langle \theta, \beta_{-1-i}\CY^3(\gamma,y)\alpha_jw\rangle y^{-1-j+i}.  
\end{array}$$
Furthermore, from  
$$\begin{array}{l}
\langle \theta, \CY^3(\alpha_{[-1]}\beta_{[-1]}[v_su],y)w\rangle \cr
\mbox{}\quad=\sum_{i,j\in \bN}\langle \theta, \alpha_{-1-j}\beta_{-1-i}\CY^3(\gamma,y)w\rangle y^{i+j}
+\sum_{i,j\in \bN}\langle \theta, \alpha_{-1-j}\CY^3(\gamma,y)\beta_iw \rangle y^{j-i-1} \cr
\mbox{}\qquad +\sum_{i,j\in \bN}\langle \theta, \beta_{-1-i}\CY^3(\gamma,y)\alpha_jw\rangle y^{-j-1+i}
+\sum_{i,j\in \bN}\langle \theta, \CY^3(\gamma,y)\beta_j\alpha_i w\rangle y^{-2-i-j}, 
\end{array}$$
we get 
$$\begin{array}{l}
{\rm LHS}\mbox{ of (4.10)}-\langle \theta, \CY^3(\alpha_{[-1]}\beta_{[-1]}[v_su],y)w\rangle \cr
=P(\sum_{i,j\in \bN}\langle \theta, \beta_{-2-i-j}[\alpha_i,\CY^3(\gamma,y)]w\rangle y^j)
-P(\sum_{i,j\in \bN}\langle \theta, [\alpha_{-1-i},\CY^3(\gamma,y)]\beta_{i+j}w\rangle y^{-1-j}) \cr
=P(\sum_{i,j,k\in \bN}\binom{i}{k}\langle \theta, \beta_{-2-i-j}\CY^3(\alpha_{[k]}\gamma,y)w\rangle y^{j+i-k}) \cr
\mbox{}\qquad -P(\sum_{i,j,k\in \bN}\binom{-1-i}{k}\langle \theta, \CY^3(\alpha_{[k]}\gamma,y)\beta_{i+j}w\rangle y^{-2-j-i-k}). 
\end{array}$$
On the other hand, we obtain  
$$\begin{array}{l}
{\rm RHS}\mbox{ of (4.10)}-\langle \theta, \CY^3(\alpha_{[-1]}\beta_{[-1]}\gamma,y)w\rangle
=P(\sum_{i\in \bN} \langle \theta, \CY^3(\alpha_{[-2-i]}\beta_{[i]}\gamma,y)w)\cr
\mbox{}\quad =P(\sum_{i,j\in \bN}\binom{-2-i}{j}(-1)^j\langle \theta, \alpha_{-2-i-j}\CY^3(\beta_{[i]}\gamma,y)w\rangle y^j)\cr
\mbox{}\qquad -P(\sum_{i,j\in \bN}\binom{-2-i}{j}(-1)^{j-2-i}\langle \theta, 
\CY^3(\beta_{[i]}\gamma,y)\alpha_jw\rangle y^{-2-i-j}).
\end{array}$$
which coincides with the above, since we obtain  
$$\begin{array}{l}
P(\sum_{i,j\in \bN}\binom{-2-i}{j}(-1)^j\langle \theta, \alpha_{-2-i-j}\CY^3(\beta_{[i]}\gamma,y)w\rangle y^j)\cr
\mbox{}\qquad =
P(\sum_{i,j,k\in \bN}\binom{i}{k}\langle \theta, \beta_{-2-i-j}\CY^3(\alpha_{[k]}\gamma,y)w\rangle y^{j+i-k}) \end{array}$$
from 
$$\begin{array}{l}
\sum_{i,j,k\in \bN}\binom{i}{k}X^{-2-i-j}Y^k Z^{j+i-k}=\sum_{i,j\in \bN}X^{-2-i-j}Z^j(Z+Y)^i \cr
\mbox{}\quad=X^{-1-i}(Z+Y)^i\iota_{X,Z}\{(X-Z)^{-1}\} 
=\iota_{X,Z,Y}\{(X-Z-Y)^{-1}(X-Z)^{-1}\}\cr
\mbox{}\quad =\sum_{i\in \bN} \iota_{X,Z}\{(X-Z)^{-2-i}\}Y^i=
\sum_{i,j\in \bN}\binom{-2-i}{j}(-1)^jX^{-2-i-j}Y^iZ^j, \end{array}$$
where we denote $\langle \theta, \alpha_a\CY^3(\beta_{[b]}\gamma,y)w\rangle y^c$ 
by $X^aY^bZ^c$ and $b\in \bN$ and 
we also have 
$$\begin{array}{l}
P(\sum_{i,j\in \bN}\binom{-2-i}{j}(-1)^{j-2-i}\langle \theta, 
\CY^3(\beta_{[i]}\gamma,y)\alpha_jw\rangle y^{-2-i-j})\cr
\mbox{}\qquad =P(\sum_{i,j\in \bN}\binom{-2-i}{j}(-1)^{j-2-i}\langle \theta, 
\CY^3(\beta_{[i]}\gamma,y)\alpha_jw\rangle y^{-2-i-j})\end{array}$$ 
from  
$$\begin{array}{l}
\sum_{i,j,k\in \bN}
\binom{-1-i}{k}X^kY^{i+j}Z^{-2-i-j-k}=\sum_{i\in \bN}(Z+X)^{-1-i}\sum_{j\in \bN} Y^iZ^{-1-j}Y^j \cr
\mbox{}\quad=\iota_{Z-Y,X}\{(Z-Y+X)^{-1}\}\iota_{Z,Y}\{(Z-Y)^{-1}\} 
=\sum_{i\in \bN} (-1)^i \iota_{Z,Y} 
\{(Z-Y)^{-2-i}\}X^i \cr
\mbox{}\quad =\sum_{i,j\in \bN}\binom{-2-i}{j}(-1)^{-2-i+j}X^iY^jZ^{-2-i+j},
\end{array}$$
where we replace $\langle \theta, \CY^3(\beta_{[a]}\gamma,y)\alpha_bw\rangle y^c$ 
with $X^aY^bZ^c$ and $a,b\in \bN$. 

This completes the proof of Theorem \ref{Vmod}.
\prend

We now start the proof of Theorem \ref{Iso}. 
Since $A\tilde{\boxtimes} B$ is an $\bN$-gradable modules, 
by Theorem \ref{Logf}, $\CY^{3}$ and $\CY^{4}$ are branches of some logarithmic intertwining 
operators of type $\binom{A(BC)}{A\tilde{\boxtimes}B\,\, C}$ and $\binom{A\tilde{\boxtimes}B}{A\,\, B}$ on $\CD^2$. 
By the construction of $\CY^{3},\CY^{4}$, 
for any $\theta\in (A(BC))^{\vee}$, there are 
$v\in A, u\in B, w\in C$ such that 
$\langle \theta, \CY^{3}(\CY^{4}(v,x-y)u,y)w\rangle\not=0$, 
that is, the set of coefficients of 
$(x-y)^ry^s\log^h(x-y)\log^k(y)$ in 
$\{\CY^{3}(\CY^{4}(v,x-y)u,y)w\mid v\in A,u\in B, w\in C\}$ spans $A(BC)$. 
From the universal properties of fusion products $(A\boxtimes B,\CY^{AB})$ and $((A\tilde{\boxtimes}B)\boxtimes C,\CY^{(A\tilde{\boxtimes}B)\boxtimes C})$ and 
$((A\boxtimes B)\boxtimes C, \CY^{(AB)C})$,  
there is a surjective $V$-homomorphism 
$\phi_{[AB]C}: (A\boxtimes B)\boxtimes C \to A\boxtimes (B\boxtimes C)$ such that 
$$\begin{array}{rl}
\langle \theta, \phi_{[AB]C}(\CY^{(AB)C}(\CY^{AB}(v,x-y)u,y)w\rangle=&\langle \theta, \CY^{3}(\CY^{4}(v,x-y)u,y)w\rangle \cr 
=&\langle \theta, \CY^{A(BC)}(v,x)\CY^{BC}(u,y)w\rangle
\end{array}$$
on $\CD^2$ 
for any $v\in A, u\in B, w\in C$ and $\theta\in (A(BC))^{\vee}$. 
Similarly, by starting from $\langle \theta',  \CY^{(AB)C}(\CY^{(AB)}(v,x-y)u,y)w\rangle$, we also have a surjective homomorphism $\phi_{A[BC]}: (A\boxtimes B)\boxtimes C \to A\boxtimes (B\boxtimes C)$ satisfying  
$$\langle \theta', \phi_{A[BC]}(\CY^{A(BC)}(v,x)\CY^{BC}(u,y)w)\rangle
=\langle \theta', \CY^{(AB)C}(\CY^{AB}(v,x-y)u,y)w\rangle $$
on $\CD^2$ for all $v,u,w,\theta'$. 
From the construction, it is easy to see that  
$\phi_{A[BC]}\circ\phi_{[AB]C}=1$ and $\phi_{[AB]C}\circ\phi_{A[BC]}=1$. 

This completes the proof of the associativity law of fusion products. \qed

\begin{cry} 
If $V^{\vee}\cong V$ and $A$ and $A^{\vee}$ are both $C_1$-cofinite 
$\bN$-gradable modules, 
 then 
for any non-zero $C_1$-cofinite $\bN$-gradable $V$-module $B$, 
$A\boxtimes B$ and $B\boxtimes A$ are not zero. 
\end{cry}

\pr There is a surjective intertwining operator 
of type $\binom{V^{\vee}}{A\, \, A^{\vee}}$ which comes from the $V$-module structure on $A$ by skew-symmetry and duality. 
Hence $B\cong B\boxtimes V$ is a homomorphism image of 
$B\boxtimes (A\boxtimes A^{\vee})\cong (B\boxtimes A)\boxtimes A^{\vee}$. 
Therefore, $(B\boxtimes A)\boxtimes A^{\vee}\not=0$ and $B\boxtimes A\not=0$. Similarly, we get $A\boxtimes B\not=0$.

\section{Pentagon axiom}
\begin{thm}[Pentagon axiom] 
 For $C_1$-cofinite $\bN$-gradable modules $A, B, C, D$, 
by using the above isomorphisms $\phi_{[\ast] \ast}$, 
we have the following commutative diagram. 
$$\begin{array}{c}
((A\boxtimes B)\boxtimes C)\boxtimes D     \xrightarrow{\phi_{([AB]C)}\times 1_D}                                                                 
 (A\boxtimes (B\boxtimes C))\boxtimes D 
\xrightarrow{\phi_{[A(BC)]D}} 
 A\boxtimes ((B\boxtimes C)\boxtimes D) \cr
\mbox{}\hspace{2.5cm} \searrow{}^{\phi_{[(AB)C]D}}   
\mbox{}\qquad    \hspace{4.5cm}      
  \swarrow{{}_{1_A\times \phi_{[BC]D}}}         \cr
 (A\boxtimes B)\boxtimes (C\boxtimes D) 
\xrightarrow{\phi_{[AB](CD)} }
 A\boxtimes (B\boxtimes (C\boxtimes D))  
\end{array}$$
\end{thm}

\pr

We consider a domain 
$\CD^3\!=\!\{(x,y,z)\!\in\! \bC^3 \mid |x|\!>\!|y|\!>\!|z|\!>\!|x\!-\!z|\!>\!|y\!-\!z|\!>\!|x\!-\!y|\!>\!0, \mbox{ and }
x,y,z,x\!-\!y,y\!-\!z,x\!-\!z\not\in \bR^{\leq 0}\}$. We note $(7,6,4)\in \CD^3\not=\emptyset$.  
Simplify the notation, we omit the notation $\boxtimes$.
Let $\theta\in (A(B(CD)))^{\vee}, v\in A, u\in B, w\in C, d\in D$ 
and let $\vec{\xi}$ denote a quadruple $(v,u,w,d)\in A\times B\times C\times D$.  
We fix all intertwining operators for fusion products, say  
$\CY^{AB},\CY^{BC},\CY^{CD},\CY^{B(CD)},\CY^{(BC)D},\CY^{A(B(CD))},\ldots$, and 
consider their five-point correlation functions:
$$\begin{array}{l}
F^{A(B(CD))}(\theta^1, \vec{\xi};x,y,z)
=\langle \theta^1,\CY^{A(B(CD))}(v,x)\CY^{B(CD)}(u,y)\CY^{CD}(w,z)d\rangle, \cr 
F^{A((BC)D)}(\theta^2, \vec{\xi};x,y\!-\!z,z)
=\langle \theta^2,\CY^{A((BC)D)}(v,x)\CY^{(BC)D}(\CY^{BC}(u,y\!-\!z)c,z)d\rangle, \cr
F^{(A(BC))D}(\theta^3, \vec{\xi};x\!-\!z,y\!-\!z,z)=\langle \theta^3, \CY^{(A(BC))D}(\CY^{(A(BC)}(v,x\!-\!z)\CY^{BC}(u,y-z)w,z)d\rangle,\cr
F^{((AB)C)D}(\theta^4,\vec{\xi};x\!-\!y,y\!-\!z,z)\!=\!\langle \theta^4,\CY^{((AB)C)D}(\CY^{(AB)C}(\CY^{AB}(v,x\!-\!y)u,y\!-\!z)w,z)d\rangle, \!\mbox{ and}\cr
F^{(AB)(CD)}(\theta^5, \vec{\xi};x\!-\!y,y,z)=\langle \theta^5,\CY^{(AB)(CD)}(\CY^{AB}(v,x\!-\!y)u, y)\CY^{CD}(w,z)d\rangle  
\end{array}\eqno{(5.1)}$$ 
for $\theta^1\in (A(B(CD)))^{\vee}, \theta^2\in (A((BC)D))^{\vee}, \theta^3\in 
((A(BC))D)^{\vee}, \theta^4\in (((AB)C)D)^{\vee}, \theta^5\in 
((AB)(CD))^{\vee}$. 
By the same way as we did for four-point correlation functions, it is not difficult to see that 
these five-point correlation functions are all locally normal convergent on 
$\{(x,y,z)\!\in\! \bC^3 \mid |x|\!>\!|y|\!>\!|z|\!>\!|x\!-\!z|\!>\!|y\!-\!z|\!>\!|x\!-\!y|\!>\!0\}$. 
We next choose their principle branches  
$$\widetilde{F}^{A(B(CD))}(\theta^1,\vec{\xi};x,y,z), 
\widetilde{F}^{A((BC)D)}(\theta^2,\vec{\xi};x,y-z,z), \ldots, 
\widetilde{F}^{(AB)(CD)}(\theta^5,\vec{\xi};x-y,y,z)$$ 
on $\CD^3$ 
by taking the values of $\log(x),\log(y),\log(z),\log(x-y),\log(x-z),\log(y-z)$ which satisfy 
$-\pi<\Im(\log(x)), \Im(\log(y)), \Im(\log(z)), \Im(\log(x-y)), 
\Im(\log(x-z)), \Im(\log(y-z))<\pi$ and 
viewing $x^a, y^b, z^c, (x-y)^d, (x-z)^e, (y-z)^f$ as 
$e^{a\log(x)}$, $e^{b\log(y)}$, $e^{c\log(z)}$, $e^{d\log(x-y)}$, $e^{e\log(x-z)}$, $e^{f\log(y-z)}$, respectively. 

From our construction of isomorphisms $\phi_{[\ast,\ast]\ast}$, by Theorem 4.1, 
we have equations:
$$\begin{array}{l}
\widetilde{F}^{A(B(CD))}(\theta^1,\vec{\xi};x,y,z)
=\widetilde{F}^{A((BC)D)}((1_A\times \phi_{[BC]D}^{\ast})(\theta^1), \vec{\xi};x,y-z,z), \cr
\widetilde{F}^{A((BC)D)}(\theta^2,\vec{\xi};x,y-z,z)
=\widetilde{F}^{(A(BC))D}(\phi_{[A(BC)]D}^{\ast}(\theta^2),\vec{\xi};x-z,y-z,z), \cr
\widetilde{F}^{(A(BC))D}(\theta^3,\vec{\xi};x-z,y-z,z)
=\widetilde{F}^{((AB)C)D}((\phi_{[AB]C}^{\ast}\times 1_D)(\theta^3), \vec{\xi};x-y,y-z,z),  \cr
\widetilde{F}^{A(B(CD))}(\theta^1,\vec{\xi};x,y,z)
=\widetilde{F}^{(AB)(CD)}(\phi_{[AB](CD)}^{\ast}(\theta^1), \vec{\xi};x-y,y,z), \mbox{ and} \cr
\widetilde{F}^{(AB)(CD))}(\theta^5,\vec{\xi};x-y,y,z)
=\widetilde{F}^{((AB)C)D}(\phi_{[(AB)C]D}^{\ast}(\theta^5), \vec{\xi};x-y,y,z),
\end{array}$$
on $\CD^3$ at least, where $\phi^{\ast}:U^{\vee}\to W^{\vee}$ denotes the dual of 
$\phi:W\to U$. We hence have 
$$\begin{array}{l}
\widetilde{F}^{((AB)C)D}((\phi_{[AB]C}^{\ast}\otimes 1_D)(\phi_{[A(BC)]D}^{\ast}((1_A\otimes \phi_{[BC]D}^{\ast})(\theta^1))), \vec{\xi};x-y,y-z,z)\cr
\mbox{}\quad=\widetilde{F}^{A(B(CD))}(\theta^1,\vec{\xi};x,y,z)=\widetilde{F}^{((AB)C)D}(\phi_{[(AB)C]D}^{\ast}(\phi_{[AB](CD)}^{\ast}(\theta^1)), 
\vec{\xi};x-y,y-z,z)
\end{array}$$
for all $\theta^1\in (A(B(CD)))^{\vee}$ and $\vec{\xi}$. 
As a consequence, we obtain $1_A\otimes \phi_{[BC]D}\circ \phi_{[A(BC)]D}\circ 
\phi_{[AB]C}\otimes 1_D
=\phi_{[AB](CD)}\circ \phi_{[(AB)C]D}$, 
which implies the commutativity of the diagram.

This completes the proof of Theorem 5.1.

\end{document}